%% file: main_siam.tex
\documentclass[onefignum,onetabnum]{siamart251216}

\usepackage{amsfonts}
\usepackage{amsmath}
\usepackage{amssymb}
\usepackage{bm}
\usepackage{graphicx}
\usepackage{epstopdf}
\usepackage{caption}
\usepackage{subcaption}
\usepackage{algorithm}
\usepackage{algpseudocode}
\usepackage{setspace}
\usepackage{ulem}
\usepackage{multirow}
\usepackage{booktabs}
\ifpdf
  \DeclareGraphicsExtensions{.eps,.pdf,.png,.jpg}
\else
  \DeclareGraphicsExtensions{.eps}
\fi
   
\algnewcommand{\IIf}[1]{\State\algorithmicif\ #1\ \algorithmicthen}
\usepackage[version=4]{mhchem}

\usepackage{epsfig, algorithm, algpseudocode}
\usepackage{multirow, booktabs}
\usepackage{xr}

\title{A Reduced Order Model for First-Principles Atomistic Simulations}

\author{
  Siu Wun Cheung\thanks{Center for Applied Scientific Computing, Lawrence Livermore National Laboratory, Livermore, CA 94550 (\email{cheung26@llnl.gov}, \email{choi15@llnl.gov}, \email{oseikuffuor1@llnl.gov})} \and
  Youngsoo Choi\footnotemark[1] \and
  Jonas Kaufman\thanks{Materials Science Division, Lawrence Livermore National Laboratory, Livermore, CA 94550 (\email{kaufman22@llnl.gov})} \and
  Daniel Osei-Kuffuor\footnotemark[1] \and
  Jean-Luc Fattebert\thanks{Computational Sciences and Engineering Division, Oak Ridge National Laboratory, Oak Ridge, TN 37830 (\email{fattebertj@ornl.gov})}
}

\usepackage{amsopn}

\ifpdf
\hypersetup{
  pdftitle={ROM for atomistic simulations},
  pdfauthor={S. W. Cheung}
}
\fi

\begin{document}

\maketitle

% REQUIRED
\begin{abstract}
To leverage the redundancy between the electronic structure computed for various model parameters such as varying atomic positions,
we present a data-driven modeling framework for Kohn-Sham Density Functional Theory 
that bypasses the explicit optimization of electronic wavefunctions. 
We sample \textit{a priori} representative atomic configurations and construct a low-dimensional basis that efficiently approximates the electronic structure subspace. 
Subsequently, we employ this reduced basis 
in an iterative solver for the electronic single particle density matrix, 
thereby enabling the efficient determination of 
ground state  without iterative wavefunction optimization. 
We demonstrate the accuracy of our approach in a Born-Oppenheimer molecular dynamics of a water molecule, showing that the resulting simulations accurately reproduce key structural properties, 
such as bond lengths and bond angle, obtained from full first-principles molecular dynamics. 
This work highlights the potential of data-driven approaches to develop low-dimensional electronic structure representations for 
first-principles simulations.
\end{abstract}

% REQUIRED

% REQUIRED

\input{introduction}
\input{FOM}
\input{ROM}
\input{PinnedH2O}
\input{conclusion}
\input{appendix}

\bibliographystyle{elsarticle-num}
\bibliography{references}

\end{document}

%% file: introduction.tex
\section{Introduction}

First-principles simulations are computational methods that determine materials or molecular properties using the fundamental laws of quantum mechanics.
By solving the Schr\"{o}dinger equations or a simplified approximation of those without relying on any empirical parameters, they predict these properties from the ground up.
The electronic structure contains all the information required to determine the physical, chemical, and spectroscopic properties of an atomistic system.
Density Functional Theory (DFT) \cite{KohnSham1965,cances2023density} is often the model of choice due to its computational efficiency and its accuracy for many problems of interest.
DFT reformulates the many-body electronic problem into a system of effective single-particle equations
which can be solved for physical systems composed of over 100 atoms. 
In the Born-Oppenheimer approximation, atomic nuclei are treated as classical particles while the electronic structure around them is modeled at the quantum level.
Representing the quantum electrons accurately can require a significant number of parameters.
Using Linear Combination of Atomic Orbitals (LCAO) is a physics-based approach where molecular orbitals (MOs) are approximated by superimposing individual atomic orbital wave functions, considered the basis set.
This enables an electronic structure representation of decent quality with a very limited number of degrees of freedom.
With such an approach however, it is difficult to systematically enrich the basis set and enhance accuracy.
Coming from the other end, the pseudo-spectral method, also known as plane-waves approach in the field, is a purely numerical discretization which can be systematically improved by increasing the basis set size \cite{ihm1979}.
Alternatives such as finite difference discretizations \cite{Briggs1995, GHOSH2017} or finite elements \cite{tsuchida1995,chen2014adaptive} approaches are also used in practice. 
For these types of discretization, the number of degrees of freedom per atom can be quite large, in the thousands.
Attempts to reduce that number of degrees of freedom while preserving accuracy have been proposed using for instance
local mesh refinement \cite{gygi1995,fattebert1999,chen2014adaptive}.
Another option proposed by Pask and Sukumar \cite{Pask2017} is to build the known atomic physics into the solution process using partition-of-unity enrichment techniques in finite element analysis.
The readers are referred to \cite{saad2010numerical} for a comprehensive overview of numerical methods for electronic structure calculations, including different discretization and eigensolver approaches.
In this paper, we look into the related question of what is the minimum number of degrees of freedom necessary to describe accurately the electronic structure of a molecule.

Molecular dynamics (MD) \cite{frenkel2023understanding} is a crucial tool in modern science for understanding the behavior of atoms and molecules over time, providing insights into material properties, chemical reactions, and biological processes. However, these simulations often rely on empirical potentials, which limit their predictive power and transferability.
The core computational task of first-principles simulations lies in the need to solve a complex, 
nonlinear electronic structure eigenvalue KS problem. 
Solving the KS equations requires first their discretization in a finite basis set or on a mesh.
While Linear Combinations of Localized Orbitals (LCAO) is a popular approach in the chemistry community, here we focus primarily on systematically improvable discretizations.
Among those the Plane Waves (PW) approach, a pseudo-spectral approach, is certainly the most used in the condensed matter community.
Here we focus on real-space finite difference discretizations, an alternative approach which allows more flexible boundary conditions and provides some advantages for parallel distributed implementations \cite{Briggs1995, GHOSH2017}.
Solving the discretized KS equations typically relies on iterative solvers given the sparse and large scale nature of the matrices involved.
The non-linearity of the KS operator is often dealt with by updating the Hamiltonian in an outer self-consistent (SC) loop, and updating the solution using an iterative solver in an inner loop for a frozen linear Hamiltonian.
Various solutions have been proposed to solve approximately the eigenvalue problem in the inner loop, including the Davidson method \cite{davidson1975iterative}, the Residual Minimization Method with Direct Inversion in the Iterative Subspace (RMM-DIIS) \cite{Kresse1996}, the Conjugate Gradient \cite{payne1992iterative}, the Chebyshev-filtered subspace iteration \cite{Zhou2006}, and the orbital minimization methods
\cite{lu2017preconditioning}. 
Alternative formulations based on gradient flows have also been
developed for the Kohn--Sham ground-state problem
\cite{dai2020gradient}.
In addition, and specially relevant to the methodology described in this paper, {\it direct} solvers which treat and solve the KS equations directly with a nonlinear solver have been proposed \cite{marzari1997ensemble, FATTEBERT2010abpg}.

In this paper, we explore a data-driven ROM approach at the electronic structure level.
We develop and evaluate a technique to reuse the electronic structure precomputed at a series of sampling points and facilitate the computation of the electronic structure for previously unseen atomic configurations, with the application of projection-based reduced-order models (ROMs). 
These techniques provide low-dimensional approximations of large-scale systems by utilizing low-dimensional representations of the solution manifold, thus significantly reducing the dimensionality of the problem.
Such low-dimensional representations can be constructed by compressing high-fidelity snapshot data through linear techniques like proper orthogonal decomposition (POD) \cite{berkooz1993proper}, balanced truncation \cite{safonov1989schur}, and the reduced basis method \cite{rozza2008reduced}.
Detailed surveys on classical projection-based ROM techniques are available in \cite{gugercin2004survey, benner2015survey}.
While projection-based ROMs are widely used for time-dependent nonlinear problems across various scientific and engineering fields, our particular focus here is on the distinct challenge of electronic structure calculations. 
A foundational contribution in this area is the work of Canc\'{e}s et al. \cite{cances2007feasibility}, which demonstrated the feasibility and competitiveness of a reduced basis approach for rapid electronic structure calculations in quantum chemistry. This work proved that the solution manifold of the nonlinear Kohn-Sham equations could be efficiently captured by a low-dimensional RB space.
Projection-based approaches have also been developed for Kohn--Sham DFT in the context of embedding, where the electronic structure problem is formulated in terms of projected subspaces \cite{lin2019projection}.

This paper introduces a novel data-driven ROM approach 
that specifically addresses the electronic structure problem in first-principles simulations. 
We focus on wavefunction-based discretizations of the KS equations, that is numerical approaches which rely on representing the electronic wavefunctions in a systematically improvable numerical basis such plane waves (pseudo-spectral) or finite elements, or rely on representing these wavefunctions on a real-space mesh (finite differences), as opposed to using Linear Combination of Atomic Orbitals (LCAO).
Unlike traditional approaches that explicitly optimize electronic wavefunctions in a high-dimensional search space, our method leverages an offline-computed, low-dimensional basis derived from snapshots of representative electronic wavefunctions. 
In the subsequent online stage, we employ this reduced basis to develop an iterative solver for the single particle density matrix in that basis, thereby enabling the evaluation of electronic structure dependent physical observables
without iterative wavefunction optimization. 
Each iteration involves solving an eigenvalue problem projected to a reduced subspace as studied in \cite{cheung2026theory}. 
From the solution of that iterative solver, physical observables can be evaluated.
In particular, atomic forces can be computed and used to propagate the atoms treated as classical particles for First-Principles Molecular Dynamics (FPMD) in the Born-Oppenheimer approximation. 
Atomic forces accuracy and energy conservation in FPMD are good fidelity tests for the converged electronic structure and will be used in this paper.
We demonstrate the accuracy of this ROM strategy through a detailed example of a water molecule, showcasing its ability to reproduce molecular trajectories and structural properties obtained from full first-principles simulations. 

The rest of the paper is organized as follows. 
In Section~\ref{sec:fom}, we describe the FPMD approach, including the electronic structure calculation by KS DFT and the ionic dynamics. 
Next, we introduce the reduced order modeling molecular dynamics (ROM-MD) framework which utilizes a ROM subspace approximation for KS DFT in Section~\ref{sec:rom}. 
We provide numerical results on an empirical example of pinned water molecular using the ROM approach in Section~\ref{sec:numerical}. 
Finally, a conclusion is given. 

%% file: FOM.tex
\section{First-Principles Atomistic Simulations}
\label{sec:fom}
 
In this section, we outline the fundamental framework of first-principles simulations. 
We begin by detailing the core electronic structure theory used to determine the ground state energy and electronic density. 
Following this, we describe how atomic forces are derived and can be used in molecular dynamics to propagate the ions through time as classical particles in the Born-Oppenheimer approximation.

\subsection{Electronic structure calculation}
\label{sec:electronic-structure}

KS DFT \cite{KohnSham1965} is a widely used first-principles model in electronic structure calculations. 
In practice, 
the electronic ground state of a physical system can be described by a system of 
one-particle electronic wavefunctions that minimize the KS total energy. 
In this work, to simplify the treatment of the core electrons, which are assumed to be tightly bound and chemically inert, we adopt the frozen core approximation. 
This approach replaces the singular Coulomb potential of the nuclei and the core electrons with smooth pseudopotentials that act on the valence electrons. 
Consequently, the electronic structure calculation problem is effectively reduced to the computation of the valence electrons. 
Second, we consider closed-shell systems, 
and assume
each spatial orbital $\phi_i$ can accommodate two valence electrons of opposite spin, spin up and spin down. 
As a result, the number of occupied orbitals $N_0$ is half the number of valence electrons.

Let $N$ be the dimension of the search trial subspace of the electronic wavefunctions, 
where $N_0 \leq N$, which is spanned by $N$ 
linearly independent real-valued wavefunctions. 
Here we adopt a general formulation that expresses the solution of the KS equations as a general set of wavefunctions that are not necessarily eigenfunctions of the KS Hamiltonian\cite{fattebert2000}.
The electronic density is defined as
\begin{equation}
\rho(\mathbf{r}) = 2 \sum_{i,j=1}^{N} X_{ij} \phi_i(\mathbf{r}) \phi_j(\mathbf{r}),
\label{eq:electronic_density}
\end{equation}
where $\mathbf{X} = (X_{ij})_{i,j=1}^N \in \text{Sym}_N(\mathbb{R})$ is the single particle density matrix (DM) which describes the occupation of the electronic wavefunctions, subject to to constraints $Tr(\mathbf{X}) = N_0$, 
and all its eigenvalues are in the interval $[0,1]$. 
Note that in the case where the electronic wavefunctions $\phi_i$ are eigenfunctions of the DFT Hamiltonian, 
$X$ is a diagonal matrix with the occupation numbers $\{ f_i \}_{i=1}^N$ as diagonal elements.
The motivation for using a general non-diagonal matrix $X$ will become clearer later on in Section \ref{sec:rom-electronic-structure}

For a molecule composed of $N_I$ atoms located in positions $\{\mathbf{R}_I\}_{I=1}^{N_I}$ in a computational domain $\Omega$, 
the KS energy functional $E_{KS}$, depends on the wavefunctions $\{\phi_i\}_{i=1}^N$ and the density matrix $\mathbf{X}$, is given by (in atomic units) 
\begin{equation}
\begin{split}
E_{KS}[\{\phi_i\}_{i=1}^N, \mathbf{X}; \{\mathbf{R}_I\}_{I=1}^{N_I}] 
& = T_s[\{\phi_i\}_{i=1}^N, \mathbf{X}] + E_H[\{\phi_i\}_{i=1}^N, \mathbf{X}] + \\ 
& \quad \quad E_{\text{xc}}[\{\phi_i\}_{i=1}^N, \mathbf{X}] + E_{ps}[\{\phi_i\}_{i=1}^N, \mathbf{X}; \{\mathbf{R}_I\}_{I=1}^{N_I}].
\end{split}
\label{eq:Kohn-Sham-energy}
\end{equation}
The first term is the kinetic energy of the non-interacting electrons, defined as
\begin{equation}
T_s[\{\phi_i\}_{i=1}^N, \mathbf{X}] 
= 2
\sum_{i,j=1}^N X_{ij} \int_\Omega \phi_i(\mathbf{r}) \left( -\frac{1}{2} \nabla^2 \right) \phi_j(\mathbf{r}) d\mathbf{r}.
\end{equation}
The second term is the Hartree energy, defined as
\begin{equation}
E_H[\{\phi_i\}_{i=1}^N, \mathbf{X}] = 
\frac{1}{2} \int_\Omega \int_\Omega \frac{\rho(\mathbf{r}) \rho(\mathbf{r'})}{|\mathbf{r} - \mathbf{r'}|} d\mathbf{r} d\mathbf{r'}, 
\end{equation}
which accounts for the classical electrostatic repulsion between the electrons.
The third term is the exchange-correlation energy, defined as 
\begin{equation}
E_{xc}[\{\phi_i\}_{i=1}^N, \mathbf{X}] = 
2 \sum_{i,j=1}^N X_{ij} \int_\Omega \phi_i(\mathbf{r}) V_{xc}(\mathbf{r}) \phi_j(\mathbf{r}) d\mathbf{r},
\end{equation}
which accounts for the many-body effects of electron-electron interactions that are not captured by the simple electrostatic (Hartree) energy in the non-interacting Kohn-Sham system. It includes the effects of the Pauli Exclusion Principle (exchange) and the dynamic correlation between electrons due to their instantaneous Coulomb repulsion, as well as a correction to the kinetic energy of the non-interacting system. Common exchange-correlation models include the local density approximation (LDA) and the generalized gradient approximation (GGA). Among GGA functionals, the Perdew–Burke–Ernzerhof (PBE) formulation \cite{perdew1996generalized} is widely used due to its favorable accuracy and compatibility with grid-based methods.
The last term is the energy of the electrons due to the potential created by the ionic cores located at positions $\{\mathbf{R}_I\}_{I=1}^{N_I}$, defined as
\begin{equation}
E_{ps}[\{\phi_i\}_{i=1}^N, \mathbf{X}; \{\mathbf{R}_I\}_{I=1}^{N_I}] = 
2 \sum_{i,j=1}^N X_{ij} \int_\Omega \phi_i(\mathbf{r}) V_{ps}(\mathbf{r}; \{\mathbf{R}_I\}_{I=1}^{N_I}) \phi_j(\mathbf{r}) d\mathbf{r},
\end{equation}
which represents the attractive electrostatic interaction between the valence electrons and the positively charged ion cores. 
These ionic cores include the nucleus and core electrons, but in a smoother form that simplifies the treatment of the chemically active valence electrons. 

The Euler-Lagrange equations associated with the minimization of $E_{KS}$ under an orthornomation constraint for the wavefunctions lead to an eigenvalue problem
\begin{equation}
    H_{KS}\psi_i=\epsilon_i\psi_i, \; i=1,\dots,N
\end{equation}
Not that each eigenfunctions $\psi_i$ can be expressed as a linear combination of the general function $\phi_i, i=1,\dots,N$ introduced earlier.
The Hamiltonian operator $H$ takes the form 
\begin{equation}
H_{KS}[\rho; \{\mathbf{R}_I\}_{I=1}^{N_I}] = 
-\frac{1}{2}\nabla^2 + v_H[\rho] + \mu_{xc}[\rho] + V_{ps}[\{\mathbf{R}_I\}_{I=1}^{N_I}].
\label{eq:Hamiltonian}
\end{equation}
Here $v_H$ is the Hartree potential which represents the Coulomb potential due to the electronic charge density $\rho$, and $\mu_{xc} = \delta E_{xc} / \delta \rho$ is the exchange-correlation potential. 
We remark that the Hartree potential is part of the the total electrostatic potential
and can be computed by solving a Poisson equation.
We can define a neutral total charge density 
the sum of the electronic density $\rho$
and the core charge density $\rho_s$
, which is defined as the sum over all ions
\cite{FATTEBERT2003}
\begin{equation}
\rho_s(\mathbf{r}) = \sum_{I=1}^{N_I} \dfrac{-Z_I}{(\sqrt{\pi} r_c^I)^3} \exp\left(-\dfrac{\vert \mathbf{r} - \mathbf{R}_I \vert^2}{(r_c^I)^2}\right).
\end{equation}
The Coulomb potential due to the ionic core charge density is given by 
\begin{equation}
v_s(\mathbf{r}) = \sum_{I=1}^{N_I} \dfrac{-Z_I}{\vert \mathbf{r} - \mathbf{R}_I \vert} \text{erf}\left(-\dfrac{\vert \mathbf{r} - \mathbf{R}_I \vert^2}{r_c^I}\right).
\end{equation}
Together, $v_H + v_s$ is the solution of the Poisson equation with total neutral source charge $\rho + \rho_s$
\begin{equation}
-\nabla^2 (v_H(\mathbf{r}) + v_s(\mathbf{r})) = 4\pi (\rho(\mathbf{r}) + \rho_s(\mathbf{r})).
\end{equation}

\subsection{Electronic structure solver}

We consider a finite difference discretizations of the wavefunction with $M$ spatial degrees of freedom, 
which defines a mapping $\mathcal{J}$ from a wavefunction matrix $\boldsymbol{\Phi} \in \mathbb{R}^{M \times N}$ to 
a set of discretized trial wavefunctions on the grid $\{ \phi_i \}_{i=1}^N = \mathcal{J}(\boldsymbol{\Phi})$, and defines 
the discrete Kohn--Sham energy functional $E_h$ by 
\begin{equation}
E_h[\boldsymbol{\Phi}, \mathbf{X}] := E_{KS}[\mathcal{J}(\boldsymbol{\Phi}), \mathbf{X}].
\label{eq:discrete-energy}
\end{equation}
The discretization also results in the discrete Hamiltonian matrix by $\mathbf{H}(\rho) \in \text{Sym}_M(\mathbb{R})$. 
A direct solver for the electronic structure calculation would involve two loops, where the outer loop updates the wavefunctions $\boldsymbol{\Phi}$, and the inner loop solves for the density matrix $\mathbf{X}$ \cite{marzari1997ensemble,YANG2006709, fattebert2022robust}.
However, when solving a problem with a wide gap between the eigenvalue associated with the highest occupied state and the eigenvalue associated with the lowest unoccupied state, simpler and computationally cheaper solvers can be used that involve only the fully occupied wavefunctions.
In that case $N=N_0$ and we can use a 
iterative optimization procedure with gradient-based corrections
\begin{equation}
    \mathbf{D} = 
    \mathbf{K} (\mathbf{H} \mathbf{\Phi} - \mathbf{\Phi} (\mathbf{\Phi}^\top \mathbf{H} \mathbf{\Phi})),
    \label{eq:gradient}
\end{equation}
where $\mathbf{K} \in \mathbb{R}^{M \times M}$ is a given preconditioner. 
Note that in practice, neither the Hamiltonian $\mathbf{H}$ nor the preconditioner $\mathbf{K}$ will be stored as matrices.
Instead, only the action of these discrete operators on vectors will be implemented.
We remark that, as an alternative to updating the wavefunctions in an additive fashion, one could
extend the search subspace spanned by $[\mathbf{\Phi}, \mathbf{D}]$ 
and optimize $\mathbf{X}$ in that subspace \cite{YANG2006709, fattebert2022robust}. 
For the FOM, we will use the Accelerated Block Preconditioned Gradient (ABPG)
proposed in \cite{FATTEBERT2010abpg}.
This algorithm is formulated for a set of non-orthogonal wavefunctions to facilitate the use of the Anderson extrapolation scheme.
In that case, the Gram matrix $\mathbf{S}=\mathbf{\Phi}^\top \mathbf{\Phi}$ is to be computed, 
and its inverse inserted into our mathematical formulation in the role of the DM $\mathbf{X}$, 
while the electronic density can be written as 
\begin{equation}
    \rho(\mathbf{r}) = 2 \sum_{i,j=1}^{N_0} (\mathbf{S}^{-1})_{ij} \phi_{i}(\mathbf{r})\phi_{j}(\mathbf{r}).
\label{eq:density}
\end{equation}
The preconditioned gradient of the KS functional becomes
\begin{equation}
    \mathbf{P} = 
    \mathbf{K} (\mathbf{H} \mathbf{\Phi} - \mathbf{\Phi} \mathbf{S^{-1}}(\mathbf{\Phi}^\top \mathbf{H} \mathbf{\Phi})),
    \label{eq:gradientS}
\end{equation}
A simplified version of the ABPG algorithm, where the Anderson extrapolation uses only one previous step to improve convergence, 
is summarized in Algorithm~\ref{alg:es-solver}. 
A two-levels multigrid V-cycle is used as preconditioner \cite{FATTEBERT2003}.
The readers are referred to \cite{FATTEBERT2010abpg} for details and handling of special cases.

\begin{algorithm}
\caption{Gradient-based optimization for electronic structure calculation}
\begin{algorithmic}[1]
\State \textbf{Input:} 
maximum number of iterations $K_{\text{ES}}$, convergence tolerance $\delta_{\text{ES}}$, and initial guess for wavefunction matrix $\boldsymbol{\Phi}_0 \in \mathbb{R}^{M \times N_0}$ 
\For{$k = 0, 1, \dots,$}
\If{$k = K_{\text{ES}}$}
\State \textbf{break}  
\EndIf
\State Evaluate electronic density $\rho_k$ according to \eqref{eq:density} 
\State Update Hamiltonian matrix operator $\mathbf{H}_k = \mathbf{H}(\rho_k) \in \text{Sym}_{M}(\mathbb{R})$
\State Compute the overlap matrix $\mathbf{S}_k = \boldsymbol{\Phi}_k^\top \boldsymbol{\Phi}_k$
\State Compute the gradient $\mathbf{R}_k =\mathbf{H}_k \mathbf{\Phi}_k - \mathbf{\Phi}_k \mathbf{S}_k^{-1}(\mathbf{\Phi}_k^\top \mathbf{H}_k \mathbf{\Phi}_k)$
\If{$\|\mathbf{R}_k\|_2 < \delta_{\text{ES}}$}
\State \textbf{break}  
\EndIf
\State Compute the preconditioned gradient-based corrections $\mathbf{P}_k = \mathbf{K} \mathbf{R}_k$
\If{$k = 0$}
\State Update wavefunction matrix 
$\mathbf{\Phi}_{k+1} = \mathbf{\Phi}_k - \mathbf{P}_k$
\Else
\State
Compute the Anderson extrapolation coefficient
$$\theta_k = 
    \dfrac{(\mathbf{P}_k-\mathbf{P}_{k-1})^\top \mathbf{P}_k}{\|\mathbf{P}_k-\mathbf{P}_{k-1}\|_2^2}
$$
\State Update wavefunction matrix 
$\mathbf{\Phi}_{k+1} = \mathbf{\Phi}_k +\theta_k (\mathbf{\Phi}_{k-1}-\mathbf{\Phi}_k) 
+ (\mathbf{P}_k + \theta_k (\mathbf{P}_{k-1}-\mathbf{P}_k))$
\EndIf
\EndFor
\State \textbf{Output:} final wavefunction matrix $\boldsymbol{\Phi}_k \in  \mathbb{R}^{M \times N_0}$ 
and density matrix $\mathbf{X} = \mathbf{S}_k^{-1} \in \text{Sym}_N(\mathbb{R})$
\end{algorithmic}
\label{alg:es-solver}
\end{algorithm}

\subsection{Ionic Forces and Dynamics}
\label{sec:ionic-dynamics}

An accurate electronic structure calculation provides all the information needed to evaluate energetic properties of a physical system, enabling energy comparisons between different atomic configurations for instance.
Energy derivatives with respect to atomic coordinates can also be computed and provide the forces exerted on an atom from the rest of the system.
These can be used to find the geometry of lowest energy of an atomistic system.
In FPMD, the dynamics of ions are governed by these forces derived from the total energy functional. 
After the minimization of the energy functional $E_{KS}$, the forces on the ions can be evaluated and 
subsequently used to propagate ionic trajectories in time via Born-Oppenheimer molecular dynamics.
By the Hellmann–Feynman theorem, given the wavefunctions $\{\phi_i\}_{i=1}^N$ and single particle density matrix $\mathbf{X}$, the total force acting on the ion $I$ located at position $\mathbf{R}_I$ and with core charge $Z_I$ is obtained as the negative partial derivative of the free energy functional with respect to the ionic coordinates. 
The only explicit dependence of the energy on the atomic positions arises from the pseudopotential and ion-ion electrostatic terms.
This yields the force expression \cite{FATTEBERT2003}
\begin{equation}
\mathbf{F}_I = - 2\sum_{i,j=1}^N X_{ij} \int_\Omega \phi_i(\mathbf{r}) \nabla_{\mathbf{R}_I} V_{ps}(\mathbf{r}; {\mathbf{R}_I}) \phi_j(\mathbf{r})  d\mathbf{r} - \sum_{J \ne I} \frac{Z_I Z_J}{|\mathbf{R}_I - \mathbf{R}_J|^3} (\mathbf{R}_I - \mathbf{R}_J). 
\label{eq:force}
\end{equation}

The Born-Oppenheimer molecular dynamics assumes the ions evolve like classical particles surrounded by quantum electrons.
Let $\Delta t$ be the the step size, 
and $K_{\text{MD}}$ be the number of timesteps in the molecular dynamics simulation. 
Given the ionic positions $\{\mathbf{R}_I(t_k)\}_{I=1}^{N_I}$ at time $t_k = k\Delta t$, for $k = 0,1, \ldots, K_{\text{MD}}-1$, 
the electronic structure calculation problem is set up using the current ionic positions and solved for the wavefunctions $\{\phi_i\}_{i=1}^N$ as described in Algorithm~\ref{alg:es-solver}. 
For the ion $I$, the ionic force is computed using \eqref{eq:force} and its position is then updated by the Verlet algorithm 
\begin{equation}
\mathbf{R}_I(t_{k+1}) = 
\begin{cases} 
\mathbf{R}_I(t_0) + \Delta t \mathbf{V}_I(t_0) + \dfrac{\Delta t^2}{2 M_I} \mathbf{F}_I(t_0) & \text{ if } k = 0, \\
2\mathbf{R}_I(t_k) - \mathbf{R}_I(t_{k-1}) + \dfrac{\Delta t^2}{2 M_I} \mathbf{F}_I(t_k) & \text{ if } 1 \leq k \leq K_{\text{MD}} - 1,
\end{cases}
\label{eq:verlet}
\end{equation}
where $M_I$ is its mass. 
The complete time integration scheme is summarized in Algorithm~\ref{alg:born-oppenheimer}.

\begin{algorithm}
\caption{Born-Oppenheimer molecular dynamics}
\begin{algorithmic}[1]
\State \textbf{Given:} initial ionic positions $\{\mathbf{R}_I(t_0)\}_{I=1}^{N_I}$ and velocities $\{ \mathbf{V}_I(t_0) \}_{I=1}^{N_I}$
\For{$k = 0, 1, \dots, K_{\text{MD}}-1$}
    \State Use Algorithm~\ref{alg:es-solver} to solve for the wavefunction matrix $\boldsymbol{\Phi} \in \text{St}(M, N)$ 
    \State Compute ionic forces $\{\mathbf{F}_I(t_k)\}_{I=1}^{N_I}$ using \eqref{eq:force}
    \State Evolve ionic positions $\{\mathbf{R}_I(t_{k+1})\}_{I=1}^{N_I}$ by \eqref{eq:verlet}
\EndFor
\end{algorithmic}
\label{alg:born-oppenheimer}
\end{algorithm}

%% file: ROM.tex
\section{Reduced Order Model}
\label{sec:rom}

In this section, we present the framework for our ROM approach. 
We exploit the inherent equivariance of the electronic structure with respect to rigid body motions (translations and rotations) of the atomic system  
for an effective compression and representation of the underlying wavefunctions.
We begin by describing the offline stage, where we strategically sample the relevant configuration space and 
construct a low-dimensional basis that accurately captures the electronic structure variations. 
Subsequently, we describe how electronic structure calculations are performed within this reduced basis to obtain the ground state energy and electronic density. 
Finally, we explain how the resulting forces on the ions are computed within the ROM approach and integrated over time 
enabling first-principles molecular dynamics, 
which provides a trajectory of ionic positions 
$\{\widetilde{\mathbf{R}}_I\}_{I=1}^{N_I}$ 
as an approximation of the groundtruth 
$\{\mathbf{R}_I\}_{I=1}^{N_I}$.

\subsection{Sampling and subspace construction}
\label{sec:offline}

We assume that the ionic positions in the laboratory frame
$\{\widetilde{\mathbf{R}}_I\}_{I=1}^{N_I}$ 
can be transformed by a rigid body transformation $\mathcal{T}: \Omega \subset \mathbb{R}^3 \to \widehat{\Omega} \subset \mathbb{R}^3$, 
into a reference frame 
$\{\widehat{\mathbf{R}}_I\}_{I=1}^{N_I} = \{ \mathcal{T} (\widetilde{\mathbf{R}}_I)) \}_{I=1}^{N_I}$ 
which can be efficiently parametrized by a vector $\boldsymbol{\nu}$ belonging to a domain $\mathsf{D} \subseteq \mathbb{R}^{n_\nu}$.
To construct a reduced basis that accurately captures the electronic structure across this domain, we employ a sampling strategy.
We select a set of $K_{\text{train}}$ representative parameter vectors from the domain $\mathsf{D}$, 
denoted as $\mathsf{D}_{\text{train}} = \{ \boldsymbol{\nu}_i \}_{i=1}^{K_{\text{train}}}$. 
For each sampled configuration $\boldsymbol{\nu}_i$, the corresponding ionic positions 
$\{\widehat{\mathbf{R}}_I(\boldsymbol{\nu}_i)\}_{I=1}^{N_I}$ define a specific atomic system in the refernce frame. 
We then minimize the Kohn-Sham energy \eqref{eq:Kohn-Sham-energy} using the electronic structure solver in Algorithm~\ref{alg:es-solver}
to obtain the corresponding wavefunction matrix $\boldsymbol{\Phi}(\boldsymbol{\nu}_i) \in \text{St}(M, N)$. 

To construct a global reduced basis that can approximate the wavefunction for any $\boldsymbol{\nu} \in \mathsf{D}$, 
we assemble a snapshot matrix $\mathbf{Y}$ by concatenating the wavefunction matrices obtained for each training configuration:
$$
\mathbf{Y} = \left[ \boldsymbol{\Phi}(\boldsymbol{\nu}_1), \boldsymbol{\Phi}(\boldsymbol{\nu}_2), \cdots, 
\boldsymbol{\Phi}(\boldsymbol{\nu}_{K_{\text{train}}}) \right] \in \mathbb{R}^{M \times K_{\text{snapshot}}},
$$
where $K_{\text{snapshot}} = N_0 K_{\text{train}} \ll M$ is the total number of snapshot vectors. 
This snapshot matrix $\mathbf{Y}$ contains a collection of representative wavefunctions from our training set. 
To extract the most important features from these snapshots and construct a low-dimensional basis, 
we perform a singular value decomposition (SVD) of $\mathbf{Y}$, which yields 
$$\mathbf{Y} = \mathbf{U} \boldsymbol{\Sigma} \mathbf{V}^\top,$$
where $\mathbf{U} \in \text{St}(M, K_{\text{snapshot}})$ and $\mathbf{V} \in \mathbb{O}(K_{\text{snapshot}})$ are orthogonal matrices consisting of the left and right singular vectors, 
and $\boldsymbol{\Sigma}$ is a diagonal matrix containing the $K_{\text{snapshot}}$ non-negative singular values in descending order.
The left singular vectors form an orthonormal basis for the column space of $\mathbf{Y}$. 
To obtain a reduced basis of size $r \ll M$, we truncate the SVD and select the first $r$ left singular vectors corresponding to the largest singular values. 
The choice of the truncation parameter $r$ is crucial and is often guided by an energy fraction criterion. 
This criterion aims to capture a certain percentage of the total ``energy'' represented by the singular values. Specifically, we determine the smallest integer $r \geq N_0$ such that
$$
\sum_{j=1}^{r} \sigma_j^2 \geq (1 - \delta_{\text{EF}}) \sum_{j=1}^{K_{\text{snapshot}}} \sigma_j^2,
$$
where $\delta_{\text{EF}}$ is a small tolerance for the energy fraction. 
The resulting basis matrix is assembled by the first $r$ left singular vectors, i.e., 
$\mathbf{Q} = \mathbf{U} \left[ \mathbf{e}_1, \mathbf{e}_2, \cdots, \boldsymbol{\Delta}_r \right] \in \text{St}(M, r)$, 
which forms an orthonormal basis for a low-dimensional subspace that effectively approximates the space spanned by the training wavefunctions, 
and can be used to approximate the electronic wavefunctions for unseen configurations $\boldsymbol{\nu} \in \mathsf{D}$ 
without the need of gradient-based corrections. 

\subsection{Electronic structure calculation}
\label{sec:rom-electronic-structure}

Having constructed a reduced basis $\mathbf{Q} \in \text{St}(M, r)$, which represents 
the reduce basis of discretized wavefunctions $\{ q_i \}_{i=1}^r = \mathcal{J}(\mathbf{Q})$ 
that accurately spans the space of relevant electronic wavefunctions, 
we can now formulate a ROM solver for electronic structure calculations of 
the wavefunctions given new atomic configurations 
$\{\widehat{\mathbf{R}}_I(\boldsymbol{\nu})\}_{I=1}^{N_I}$ in the reference frame parametrized by $\boldsymbol{\nu} \in \mathsf{D}$. 
The core idea is to use the span of the bilinear products of the reduced basis $\{ q_i \}_{i=1}^r$ to define the electronic density 
\begin{equation}
\widehat{\rho}(\widehat{\mathbf{r}}) = 2 \sum_{i,j=1}^{r} \widehat{X}_{ij} q_i(\widehat{\mathbf{r}}) q_j(\widehat{\mathbf{r}}),
\label{eq:rom-density}
\end{equation}
with respect to the single particle density matrix (DM) 
$\widehat{\mathbf{X}} = (\widehat{X}_{ij})_{i,j=1}^r \in \mathcal{D}_r$, where the admissible set $\mathcal{D}_r$ is defined as
\begin{equation}
\begin{split}
\mathcal{D}_r 
& = \left\{ \widehat{\mathbf{X}}\in\mathrm{Sym}_r(\mathbb R) : \mathbf{0} \preccurlyeq \widehat{\mathbf{X}} \preccurlyeq \mathbf{I}, \operatorname{Tr}(\widehat{\mathbf{X}})=N_0\right\} \\
& = \left\{ \mathbf{U}\boldsymbol{\Lambda} \mathbf{U}^\top : \mathbf{U} \in \mathbb{O}(r), \boldsymbol{\Lambda} = \text{diag}(\lambda_i), \lambda_i \in [0,1] \text{ for } 1 \leq i \leq r, \operatorname{Tr}(\boldsymbol{\Lambda})=N_0 \right\},
\end{split}
\end{equation}
where is a relaxation to the set of rank-$N_0$ projectors
\begin{equation}
\widetilde{\mathcal{D}}_r 
= \left\{ \mathbf{U}\boldsymbol{\Lambda} \mathbf{U}^\top : \mathbf{U} \in \mathbb{O}(r), \boldsymbol{\Lambda} = \text{diag}(\lambda_i), \lambda_i \in \{0,1\} \text{ for } 1 \leq i \leq r, \operatorname{Tr}(\boldsymbol{\Lambda})=N_0 \right\}.
\end{equation}
The ROM basis defines the reduced Kohn--Sham energy functional $E_r$ on $\mathcal{D}_r$ by 
\begin{equation}
E_r(\widehat{\mathbf{X}}) := E_{KS}[\{q_i\}_{i=1}^{r}, \widehat{\mathbf{X}}],
\label{eq:rom-energy}
\end{equation}
As a result, the ROM electronic structure calculation problem can be recast with the density matrix as the primary variable, 
avoiding the iterative optimization of electronic wavefunctions in the full space. 

Next, we quantify the proximity of the ROM problem to the original FOM problem.
We assume the discrete Kohn--Sham energy functional $E_h$ defined in \eqref{eq:discrete-energy} 
attains global minimum $E_0$ across different dimensions $N \geq N_0$. 
Among the ground-state minimizers, we choose the orthogonal representation 
$\boldsymbol{\Psi}^\star \in \mathcal{M} = \mathrm{St}(M,N_0)$, i.e., 
\begin{equation}
E_0 = E_h[\boldsymbol{\Psi}^\star, \mathbf{I}] \leq E_h[\boldsymbol{\Phi}, \mathbf{X}], 
\text{ for all } \boldsymbol{\Phi} \in \mathbb{R}^{M \times N}, \mathbf{X} \in \mathcal{D}_{N}, N \geq N_0.
\end{equation}
On the other hand, we assume the reduced Kohn-Sham energy functional $E_r$ defined in \eqref{eq:rom-energy} attains minimum $\widehat{E}_0$ 
at $\widehat{\mathbf{X}}^\star \in \mathcal{D}_r$, i.e., 
\begin{equation}
\widehat{E}_0 = E_r(\widehat{\mathbf{X}}^\star) \leq E_r(\widehat{\mathbf{Y}}), 
\text{ for all } \widehat{\mathbf{Y}} \in \mathcal{D}_r.
\end{equation}
We show that the perturbation in the ground-state energy 
introduced by the ROM formulation can be estimated by the projection error.

\begin{theorem}
\label{thm:rom-energy}
Suppose $E_h[\cdot, \mathbf{I}] \vert_{\mathcal{M}}$ is locally $C^2$-regular at the ground-state orbitals matrix $\boldsymbol{\Psi}^\star$. 
Then there exists a constant $C>0$, independent of $r$, such that
\begin{equation}
E_0 \leq \widehat{E}_0 \leq E_0 + C\varepsilon_r^2,
\end{equation}
where $\varepsilon_r = \left\| (\mathbf I-\mathbf Q\mathbf Q^\top) \boldsymbol{\Psi}^\star \right\|_2$.
\begin{proof}
It is trivial to see that $E_0 \leq \widehat{E}_0$.
Let 
$\widetilde{\boldsymbol{\Phi}}_r = \mathbf Q\mathbf Q^\top \boldsymbol{\Psi}^\star$
denote the orthogonal projection of occupied orbitals onto the reduced space $\text{range}(\mathbf{Q})$ and 
$\widetilde{\boldsymbol{\Psi}}_r
=
\widetilde{\boldsymbol{\Phi}}_r
\left(
\widetilde{\boldsymbol{\Phi}}_r^\top
\widetilde{\boldsymbol{\Phi}}_r
\right)^{-1/2} \in \mathcal{M}$
denote its orthonormalization.
Since  
$\widetilde{\boldsymbol{\Psi}}_r \in\operatorname{span}(\mathbf Q)$,
we write $\widetilde{\boldsymbol{\Psi}}_r
= \mathbf Q\widehat{\mathbf C}_r$, where 
$\widehat{\mathbf C}_r =
\mathbf Q^\top\widetilde{\boldsymbol{\Psi}}_r$, and 
define the associated reduced density matrix
$\widetilde{\mathbf X}_r
=
\widehat{\mathbf C}_r
\widehat{\mathbf C}_r^\top \in \mathcal D_r$. 
Then we have 
\begin{equation}
E_h[\widetilde{\boldsymbol{\Psi}}_r, \mathbf{I}] = E_r(\widetilde{\mathbf X}_r) \geq E_r(\widehat{\mathbf{X}}^\star) = \widehat{E}_0. 
\end{equation}
Let $\boldsymbol{\Delta}_r = \boldsymbol{\Psi}^\star - \widetilde{\boldsymbol{\Phi}}_r$
denote the orthogonal projection of occupied orbitals onto the orthogonal complement of the reduced space. 
By direct computation, we have 
\begin{equation}
    \widetilde{\boldsymbol{\Phi}}_r^\top \widetilde{\boldsymbol{\Phi}}_r + \boldsymbol{\Delta}_r^\top \boldsymbol{\Delta}_r = \mathbf{I}.
\end{equation}
On the other hand, by definition, we have $\|\boldsymbol{\Delta}_r\|_2=\varepsilon_r$ and $\|\boldsymbol{\Delta}_r^\top \boldsymbol{\Delta}_r \|_2=\varepsilon_r^2$.
By Taylor's expansion, 
\begin{equation}
    \widetilde{\boldsymbol{\Psi}}_r 
    = (\boldsymbol{\Psi}^\star - \boldsymbol{\Delta}_r)(\mathbf{I} - \boldsymbol{\Delta}_r^\top \boldsymbol{\Delta}_r)^{-1/2} 
    = (\boldsymbol{\Psi}^\star - \boldsymbol{\Delta}_r)\left(\mathbf{I} + \dfrac{1}{2} \boldsymbol{\Delta}_r^\top \boldsymbol{\Delta}_r + O(\| \boldsymbol{\Delta}_r \|_2^4) \right),
\end{equation}
which gives
\begin{equation}
\widetilde{\boldsymbol{\Psi}}_r 
= \boldsymbol{\Psi}^\star - \boldsymbol{\Delta}_r 
+ \dfrac{1}{2} \boldsymbol{\Psi}^\star \boldsymbol{\Delta}_r^\top \boldsymbol{\Delta}_r + O(\| \boldsymbol{\Delta}_r \|_2^3), 
\end{equation}
and therefore 
$\| \widetilde{\boldsymbol{\Psi}}_r-\boldsymbol{\Psi}^\star
\|_2 \leq C\varepsilon_r$. 
Since
$\boldsymbol{\Psi}^\star$
is a ground-state solution,
it satisfies the first-order optimality condition
$\text{grad}_\mathcal{M} E_h[\boldsymbol{\Psi}^\star, \mathbf{I}] = 0$. 
Take the constant-speed minimizing geodesic $\gamma: [0,1] \to \mathcal{M}$ such that 
$\gamma(0) = \boldsymbol\Psi^\star$, 
$\gamma(1) = \widetilde{\boldsymbol\Psi}_r$
and $\| \dot{\gamma}(t) \|_F = d_\mathcal{M}(\boldsymbol\Psi^\star, \widetilde{\boldsymbol\Psi}_r)$ for $t \in (0,1)$, define 
$g(t) = E_h[\gamma(t), \mathbf{I}]$.
Then we have $g'(0) = \langle \text{grad}_\mathcal{M} E_h[\boldsymbol{\Psi}^\star, \mathbf{I}], \dot{\gamma}(0) \rangle = 0$. 
The Taylor's remainder theorem gives
\begin{equation}
g(1) = g(0) + g'(0) + \int_0^1 (1-t) g''(t) \, dt, 
\end{equation} 
which implies
\begin{equation}
E_h[\widetilde{\boldsymbol{\Psi}}_r, \mathbf{I}] - E_0 
= \int_0^1 (1-t) \text{Hess}_\mathcal{M} E_h[\gamma(t), \mathbf{I}][\dot{\gamma}(t), \dot{\gamma}(t)] \, dt.
\end{equation}
Due to the $C^2$ regularity, we have 
\begin{equation}
    \text{Hess}_\mathcal{M} E_h[\gamma(t), \mathbf{I}][\dot{\gamma}(t), \dot{\gamma}(t)] 
    \lesssim \| \dot{\gamma}(t) \|_F^2 
    \lesssim \| \dot{\gamma}(t) \|_2^2 
    = d_\mathcal{M}(\boldsymbol\Psi^\star, \widetilde{\boldsymbol\Psi}_r)^2
    \lesssim \| \boldsymbol\Psi^\star- \widetilde{\boldsymbol\Psi}_r \|_2^2
    \lesssim \varepsilon_r^2,
\end{equation}
for all $t \in (0,1)$. Therefore 
\begin{equation}
    E_h[\widetilde{\boldsymbol{\Psi}}_r, \mathbf{I}] - E_0 
    \lesssim \varepsilon_r^2. 
\end{equation}
Combining the above inequalities completes the proof.
\end{proof}
\end{theorem}
We remark that Theorem~\ref{thm:rom-energy} presents an \textit{apriori} estimates on the ground-state energies, 
but the constant $C > 0$ depends on the Riemannian smoothness of the discrete Kohn--Sham energy functional $E_h$ and is hardly computable in practice.

Another key and novel ingredient in our ROM approach is the formulation of the minimization problem which works naturally on a given subspace. 
Instead of the KS energy, we work with the Mermin free energy
functional at finite electronic temperature $T$ 
\cite{Mermin1965}. 
The reduced basis $\{ q_i \}_{i=1}^r$ is used to define the free energy objective functional solely in $\widehat{\mathbf{X}}$, i.e., 
\begin{equation}
A_r^T[\widehat{\mathbf{X}}] = E_r[\widehat{\mathbf{X}}] - T S_r[\widehat{\mathbf{X}}],
\label{eq:rom-free-energy}
\end{equation}
where $S_r$ denotes the entropy term 
\begin{equation}
S_r[\widehat{\mathbf{X}}] = -2 k_B \sum_{i=1}^r \left[f_i \ln(f_i) + (1 - f_i) \ln(1 - f_i)\right], 
\label{eq:entropy}
\end{equation}
Note that for the ground state solution that minimizes \eqref{eq:rom-free-energy}, 
the occupation numbers, $f_i, i=1,\dots,N$, follow a Fermi-Dirac thermal distribution $f_i = f_T(\epsilon_i - \mu)$, where
\begin{equation}
f_T(\epsilon) = \frac{1}{e^{\epsilon/k_B T} + 1},
\label{eq:fermi-dirac}
\end{equation}
$\epsilon_i$ is the $i$-th eigenvalue of $\mathbf{H}$,
$k_B$ is the Boltzmann constant, and 
$\mu$ is the chemical potential chosen to satisfy $\sum_i f_i = N_0$,
The rationale of working with the Mermin free energy functional is that it is a regularized formulation of the reduced KS minimization problem, and in the zero-temperature limit it recovers the same optimal occupied subspace inside the reduced basis $\{ q_i \}_{i=1}^r$, as stated in the theorem below.

\begin{theorem}
\label{thm:zero-temperature-limit}
Assume that $E_r$ is continuous on $\mathcal D_r$ and possesses a unique minimizer
\begin{equation}
\widehat{\mathbf X}^\star = \operatorname*{arg\,min}_{\widehat{\mathbf X}\in\mathcal D_r} E_r(\widehat{\mathbf X}).
\end{equation}
and the corresponding reduced density
\begin{equation}
\widehat{\rho}^\star = 2 \sum_{i,j=1}^{r} \widehat{X}^\star_{ij} q_i q_j .
\end{equation}
For every $T>0$, let
\begin{equation}
\widehat{\mathbf X}^T = \operatorname*{arg\,min}_{\widehat{\mathbf X}\in\mathcal D_r} A_r^T(\widehat{\mathbf X}).
\end{equation}
and the corresponding reduced density
\begin{equation}
\widehat{\rho}^T = 2 \sum_{i,j=1}^{r} \widehat{X}^T_{ij} q_i q_j,
\end{equation}
Then $\widehat{\mathbf X}^T \to \widehat{\mathbf X}^\star$ and 
$\widehat{\rho}^{T}\to\widehat{\rho}^\star$ in $L^2(\Omega)$ as $T \to 0$.
\begin{proof}
Let $(T_n)_{n \in \mathbb{N}}$ be a positive sequence decreasing to 0. 
Since $\mathcal D_r$ is compact, there exists a convergent subsequence of $(\widehat{\mathbf X}^{T_{n}})_{n \in \mathbb{N}}$. 
Let $(\widehat{\mathbf X}^{T_{n_k}})_{k \in \mathbb{N}}$ be a convergent subsequence and 
$\widehat{\mathbf Y}\in\mathcal D_r$. 
Since $\widehat{\mathbf X}^{T_{n_k}}$ minimizes $A_r^{T_{n_k}}$,
\begin{equation}
A_r^{T_{n_k}}(\widehat{\mathbf X}^{T_{n_k}}) \le A_r^{T_{n_k}}(\widehat{\mathbf Y}).
\end{equation}
For every occupation number $f\in[0,1]$,
\begin{equation}
0 \le -\bigl(f\log f + (1-f)\log(1-f)\bigr) \le \log2 .
\end{equation}
Hence
\begin{equation}
0 \le S_r(\widehat{\mathbf Y}) \le 2 \log 2 \, r k_B.
\end{equation}
Therefore we have 
\begin{equation}
\begin{split}
E_r(\widehat{\mathbf X}^{T_{n_k}})
\le
A_r^{T_{n_k}}(\widehat{\mathbf X}^{T_{n_k}})
\le
A_r^{T_{n_k}}(\widehat{\mathbf Y})
\le
E_r(\widehat{\mathbf Y}) + 2k_Br(\log2)\,T_{n_k} .
\end{split}
\end{equation}
Passing to the limit and using the continuity of $E_r$ yields
\begin{equation}
E_r\left(\lim_{k \to \infty} \widehat{\mathbf X}^{T_{n_k}}\right) \le E_r(\widehat{\mathbf Y}).
\end{equation}
Since $\widehat{\mathbf Y}\in\mathcal D_r$ is arbitrary, 
the limit of $(\widehat{\mathbf X}^{T_{n_k}})_{k \in \mathbb{N}}$ is a minimizer of $E_r$, 
which is $\widehat{\mathbf X}^\star$ by uniqueness.
Since the convergent subsequence was arbitrary, every convergent subsequence has limit $\widehat{\mathbf X}^\star$. 
By compactness, this implies the entire sequence $(\widehat{\mathbf X}^{T_n})_{n \in \mathbb{N}}$ converges to $\widehat{\mathbf X}^\star$.
Lastly, using the fact that the map $\widehat{\mathbf X} \mapsto 2 \sum_{i,j=1}^{r} \widehat X_{ij} q_i q_j$
is linear and continuous from $\mathrm{Sym}_r(\mathbb R)$ into $L^2(\Omega)$, 
we have $\widehat{\rho}^{T_n} \to\widehat{\rho}^\star$ as $T_n\to0$.
\end{proof}
\end{theorem}

Theorem~\ref{thm:zero-temperature-limit} establishes that the finite-temperature density matrix optimization problem solved in Algorithm~\ref{alg:dm-solver} is a consistent regularization of the reduced Kohn--Sham energy minimization problem. In particular, as the electronic temperature tends to zero, any sequence of minimizers of the reduced Mermin free energy converges to a minimizer of the reduced Kohn--Sham energy within the fixed reduced space $\mathcal V_r = \mathrm{span}\{q_1,\ldots,q_r\}$. Consequently, the finite-temperature formulation does not alter the limiting reduced ground state, but instead provides a smooth optimization landscape that is often more robust numerically.

Furthermore, suppose that the self-consistent projected Hamiltonian associated with the zero-temperature solution satisfies a spectral gap condition
\[
\lambda_{N_0} < \lambda_{N_0+1},
\]
where $\lambda_{N_0}$ and $\lambda_{N_0+1}$ denote the highest occupied and lowest unoccupied eigenvalues, respectively. In this case, the reduced Kohn--Sham minimizer is unique and is given by the rank-$N_0$ projector
\[
\widehat{\mathbf X}_0
=
\widehat{\mathbf V}_0
\begin{bmatrix}
\mathbf I_{N_0} & \mathbf 0 \\
\mathbf 0 & \mathbf 0
\end{bmatrix}
\widehat{\mathbf V}_0^\top,
\]
where the columns of $\widehat{\mathbf V}_0$ are the eigenvectors of the converged projected Hamiltonian. Therefore,
\[
\widehat{\mathbf X}_T \rightarrow \widehat{\mathbf X}_0
\qquad\text{as}\qquad T\rightarrow0,
\]
and the finite-temperature fractional occupation numbers converge to the zero-temperature occupied-state projector. This shows that the density matrix optimization procedure employed in the ROM framework recovers, in the zero-temperature limit, the same occupied subspace that would be obtained by directly solving the reduced Kohn--Sham problem. 

In \cite{marzari1997ensemble}, Marzari et al. introduced a direct optimization strategy for the Mermin free energy \eqref{eq:rom-free-energy}.
We adapt this approach by reformulating the problem to be solved within a fixed, low-dimensional search subspace spanned by the reduced basis vectors $\{q_i\}_{i=1}^{r}$.
In that case, it becomes equivalent to the Optimal Damping Algorithm (ODA) proposed by Canc\`es \cite{Cances2001}, albeit in a different basis.
Algorithm~\ref{alg:dm-solver} outlines an iterative procedure for obtaining the electronic density matrix within this reduced basis.
In practice, to reduce the computational cost, we fix the mixing parameter $\beta$. 
For problems with a large band gap like the one considered later in this article, it shows a similar convergence rate without the cost associated with the extra evaluation of the electronic density and energy needed for the line minimization.

\begin{algorithm}
\caption{Density matrix optimization in a $r$-dimensional search subspace, with $r \geq N_0$}
\begin{algorithmic}[1]
\State \textbf{Input:}
orthogonal basis matrix $\mathbf{Q} \in \text{St}(M, r)$, 
maximum number of iterations $K_{\text{DM}}$, convergence tolerance $\delta_{\text{DM}}$, 
mixing parameter $\beta$, 
and initial guess of density matrix $\widehat{\mathbf{X}}_0 \in \text{Sym}_{r}(\mathbb{R})$
\For{$k = 0, 1, \dots$}
\State Evaluate electronic density $\widehat{\rho}_k$ according to \eqref{eq:rom-density}
\State Projected Hamiltonian matrix $\widehat{\mathbf{H}}_k = \mathbf{Q}^\top \mathbf{H}(\widehat{\rho}_k) \mathbf{Q} \in \text{Sym}_{r}(\mathbb{R})$
\State Solve the symmetric eigenvalue problem $\widehat{\mathbf{H}}_k \widehat{\mathbf{V}}_k = \widehat{\mathbf{V}}_k \widehat{\boldsymbol{\mathcal{E}}}_k$
\State
Form the occupation 
matrix $\widehat{\mathbf{F}}_k = f_T(\widehat{\boldsymbol{\mathcal{E}}}_k - \nu \mathbf{I})$ computed entry-wise on the diagonal according to \eqref{eq:fermi-dirac}
\State
Update the density matrix  
$\widehat{\mathbf{X}}_{k+1} = 
(1 - \beta) \widehat{\mathbf{X}}_k
    + \beta \widehat{\mathbf{V}}_k \widehat{\mathbf{F}}_k \widehat{\mathbf{V}}_k^\top$
\If{$k+1 = K_{\text{MD}}$ or 
$|A_r(\widehat{\mathbf{X}}_{k}) - A_r(\widehat{\mathbf{X}}_{k+1})| < \delta_{\text{DM}}$} 
\State \textbf{break}  
\EndIf
\EndFor
\State \textbf{Output:} final density matrix $\widehat{\mathbf{X}}_{k+1} \in \text{Sym}_r(\mathbb{R})$
\end{algorithmic}
\label{alg:dm-solver}
\end{algorithm}

%%%%%%%%%%%%%%%%%%%%%%%%%%%%%%%%%%%%%%%%%%%%%%%%%%%%%%%%%
\subsection{Ionic Forces and Dynamics}
\label{sec:rom-ionic-dynamics}

Given the density matrix $\widehat{\mathbf{X}}$, we compute the force acting on each atom in the reference frame
using an analogue to the Hellmann–Feynman force expression \eqref{eq:force}, i.e., 
\begin{equation}
\widehat{\mathbf{F}}_I = - \sum_{i,j=1}^r \widehat{X}_{ij} \int_\Omega q_i(\widehat{\widehat{\mathbf{r}}}) \nabla_{\widehat{\mathbf{R}}_I} V_{ps}(\widehat{\mathbf{r}}; {\widehat{\mathbf{R}}_I}) q_j(\widehat{\mathbf{r}})  d\widehat{\mathbf{r}} - \sum_{J \ne I} \frac{Z_I Z_J}{|\widehat{\mathbf{R}}_I - \widehat{\mathbf{R}}_J|^3} (\widehat{\mathbf{R}}_I - \widehat{\mathbf{R}}_J).
\label{eq:rom-force}
\end{equation}

The ROM-MD simulation begins with the initial ionic positions 
set as $\widetilde{\mathbf{R}}_I(0) = \mathbf{R}_I(0)$ for all ions $1 \leq I \leq N_I$. 
Subsequently, 
given the ionic positions $\{\widetilde{\mathbf{R}}_I(t_k)\}_{I=1}^{N_I}$ at time $t_k = k\Delta t$, for $k = 0,1, \ldots, K_{\text{MD}}-1$, 
we first identify the rigid body transformation $\mathcal{T}(t_k)$ that aligns the system to the domain $\widehat{\Omega}$. 
Following this transformation, we transform the coordinate system into 
the reference frame via $\widehat{\mathbf{R}}_I(t_k) = \mathcal{T}(t_k) \widetilde{\mathbf{R}}_I(t_k)$,
set up and solve the electronic structure problem, 
obtain the density matrix $\widehat{\mathbf{X}}$ as detailed in Algorithm~\ref{alg:dm-solver}. 
For each ion $I$, the ionic force $\widehat{\mathbf{F}}_I(t_k)$ in the reference frame is computed using \eqref{eq:rom-force}, 
transformed back to the laboratory frame via 
$\widetilde{\mathbf{F}}_I(t_k) = [\mathcal{T}(t_k)]^{-1} \widehat{\mathbf{R}}_I(t_k)$,
and update the ionic positions $\{\widetilde{\mathbf{R}}_I(t_{k+1})\}$ using the same Verlet algorithm as \eqref{eq:verlet} in Section~\ref{sec:ionic-dynamics}.
The complete time integration scheme for the ROM-MD is summarized in Algorithm~\ref{alg:rom-born-oppenheimer}.

\begin{algorithm}
\caption{Reduced Born-Oppenheimer molecular dynamics}
\begin{algorithmic}[1]
\State \textbf{Given:} wavefunction basis matrix $\mathbf{Q} \in \text{St}(M,r)$ in reference frame, initial ionic positions $\{\widetilde{\mathbf{R}}_I(t_0)\}_{I=1}^{N_I}$ and velocities $\{ \widetilde{\mathbf{V}}_I(t_0) \}_{I=1}^{N_I}$ in the laboratory frame
\For{$k = 0, 1, \dots, K_{\text{MD}}-1$}
    \State Identify rigid body transformation $\mathcal{T}(t_k)$ from laboratory frame to reference frame
    \State Transform atomic coordinate into reference frame $\widehat{\mathbf{R}}_I(t_k) = \mathcal{T}(t_k) \widetilde{\mathbf{R}}_I(t_k)$ for $1 \leq I\leq N_I$
    \State Use Algorithm~\ref{alg:dm-solver} solve for the density matrix $\widehat{\mathbf{X}} \in \text{Sym}_r(\mathbb{R})$
    \State Compute ionic forces in reference frame $\{ \widehat{\mathbf{F}}_I(t_k) \}_{I=1}^{N_I}$ using \eqref{eq:rom-force}
    \State Transform ionic force into laboratory frame 
    $\widetilde{\mathbf{F}}_I(t_k) = \mathcal{T}^{-1}(t_k) \widehat{\mathbf{F}}_I(t_k)$ for $1 \leq I\leq N_I$
    \State Evolve ionic positions $\{\widetilde{\mathbf{R}}_I(t_{k+1})\}_{I=1}^{N_I}$ by \eqref{eq:verlet}
\EndFor
\end{algorithmic}
\label{alg:rom-born-oppenheimer}
\end{algorithm}

%% file: PinnedH2O.tex
\section{Molecular dynamics application: pinned water molecule}
\label{sec:numerical}

In this section, we consider an empirical example with a single water molecule, 
where the Oxygen atom is considered ``pinned'', 
meaning its position is fixed throughout the simulation. 
This can be physically imagined as the oxygen atom having an infinitely large mass 
and remaining immobile under the influence of forces. 
Consequently, we focus on the forces and dynamics of the two hydrogen atoms relative to the fixed oxygen.
The computation domain is a cube $\Omega = [-6, 6]^3$ of side length 12.0 Bohr, discretized by a $64 \times 64 \times 64$ uniform grid. 
For electronic structure calculation, we consider the valence electrons to determine its chemical properties. 
Oxygen has 6 valence electrons and 2 core electrons, while each hydrogen atom contributes 1 valence electron. 
Applying the frozen core approximation, we treat the 2 core electrons of oxygen implicitly through a pseudopotential, 
focusing solely on the valence electrons. 
Furthermore, we ignore the spin degrees of freedom of the electrons. 
This implies that the $N_0$ number of occupied orbitals in our electronic structure calculations is half the total number of valence electrons, which is $(6 + 1 + 1)/2 = 4$. 
Kohn-Sham energy is defined with the PBE exchange-correlation functional \cite{perdew1996generalized} and 
norm-conserving pseudopotentials of the SG15 family for the core electrons \cite{ONCV2013,ONCV2015}.
The electronic wavefunctions were discretized on a uniform real-space grid using a a fourth-order accuracy finite difference scheme. 
As a preconditioner $\mathbf{K}$, we use the multigrid preconditioner described in \cite{FATTEBERT2003}. 
For all the FOM solve by Algorithm~\ref{alg:es-solver}, the convergence tolerance is set to $\delta_{\text{ES}} = 10^{-8}$.

\subsection{Sampling and subspace construction}

\begin{figure}[ht!]
\centering
\includegraphics[width=0.6\linewidth]{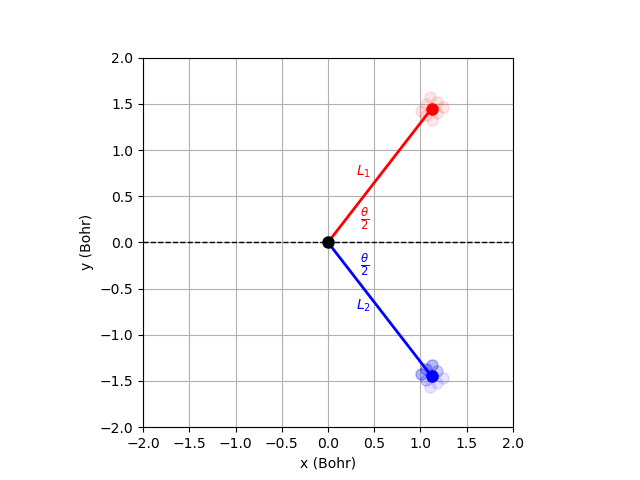}
\caption{Configurations of the water molecule parametrized by two bond lengths $(L_1, L_2)$ and one bond angle $\theta$.}
\label{fig:PinnedH2O}
\end{figure}

Next, we describe the data sampling strategy for constructing our reduced basis. 
To define a unique atomic system with rotational invariance, we assume the sampled molecule lies in the plane $z = 0$, 
the first hydrogen atom \ce{H1} lies in the first quadrant, 
the second hydrogen atom \ce{H2} lies in the fourth quadrant, 
the bond length between the oxygen \ce{O} and the first hydrogen atom \ce{H1} is not shorter than that between \ce{O} and \ce{H2}, 
and the bondangle is equally split into half by the positive $x$-axis. 
Under these assumptions, we can parametrize the molecular geometry by $\boldsymbol{\nu} = (L_1, L_2, \theta)$ using only $n_\nu = 3$ parameters: 
the bond length $L_1$ between oxygen \ce{O} and \ce{H1} , 
the bond length $L_2$ between oxygen \ce{O} and \ce{H2}, 
and the bond angle $\theta$ between the two \ce{O-H} bonds. 
We consider limited ranges for these parameters around the equilibrium geometry of 
the water molecule.
Specifically, we define the following parametric domains:
\begin{itemize}
    \item Bond length $L_1 = 1.83 s_1$ Bohr, where $s_1$ is sampled uniformly from the interval $[0.95, 1.05]$ with a step size of $\Delta_L = (1.05 - 0.95) / K_L$.
    \item Bond length $L_2 = 1.83 s_2$ Bohr, where $s_2$ is sampled uniformly from the interval $[0.95, s_1]$ with the same step size $\Delta_L$. This ensures that $L_2$ is always less than or equal to $L_1$, reflecting our initial setup where \ce{H1} has a slightly longer bond length.
    \item Bond angle $\theta = 104.5^\circ + s_\theta$, where $s_\theta$ is sampled uniformly from the interval $[-5^\circ, 5^\circ]$ with a step size of $\Delta_\theta = (5 - (-5)) / K_\theta$.
\end{itemize}
Here, $K_L$ and $K_\theta$ are the number of sampling intervals for the bond lengths and the bond angle, respectively. 
The uniform sampling in the parametric domain allows us to systematically explore a relevant portion of the configuration space. 
The number of sampled pairs $(s_1, s_2)$ is $(K_L+1) (K_L+2) / 2$, while the number of sampled angles is $K_\theta + 1$. 
Therefore, the number of unique atomic systems generated by this sampling strategy is therefore given by $K_{\text{train}} = (K_L+1) (K_L+2) (K_\theta + 1) / 2$. 
Figure~\ref{fig:PinnedH2O} depicts several configurations of the water molecule parametrized by bond lengths and bond angle in their corresponding ranges.
Since we have $N_0 = 4$ occupied wavefunctions for each atomic system, the total number of wavefunction snapshots collected for the construction of the reduced basis is 
$K_{\text{snapshot}} = 2(K_L+1) (K_L+2) (K_\theta + 1)$. 

Table~\ref{tab:rom_dimension} presents the required dimension $r$ of the reduced basis 
to achieve a specific energy fraction tolerance $\delta_{\text{EF}}$ for different sampling frequencies, 
characterized by the number of intervals $K_L$ and $K_\theta$. 
The bottom row, corresponding to $\delta_{\text{EF}} = 0$, represents the full rank of the snapshot matrix, which is equal to $K_{\text{snapshot}}$.
As expected, achieving a smaller energy fraction residual requires a larger reduced basis dimension $r$. 
Furthermore, increasing the sampling frequency generally leads to a slightly higher required reduced basis dimension for a given target accuracy tolerance $\delta_{\text{EF}}$. 
This indicates that a more comprehensive sampling of the configuration space captures a wider range of electronic states, 
thus requiring a larger basis to represent them effectively.

\begin{table}[h!]
    \centering
    \begin{tabular}{|c||cccc|}
        \hline
        \multirow{2}{*}{$\delta_{\text{EF}}$} & \multicolumn{4}{c|}{$(K_L, K_\theta)$} \\
        & $(2,2)$ & $(2,5)$ & $(5,2)$ & $(5,5)$ \\ 
        \hline
        \hline
        $10^{-2}$ & 8 & 7 & 8 & 7 \\
        $10^{-3}$ & 18 & 18 & 18 & 18 \\
        $10^{-4}$ & 34 & 35 & 36 & 36 \\
        $0$         & 72 & 144 & 252 & 504 \\
        \hline
    \end{tabular}
    \caption{Reduced basis dimension $r$ required to achieve a target energy fraction residual $\delta_{\text{EF}}$ for different sampling frequencies $K_L$ and $K_\theta$. The last row with $\delta_{\text{EF}}=0$ has the dimension $r = K_{\text{snapshot}}$.}
    \label{tab:rom_dimension}
\end{table}

\subsection{Accuracy in forces}

To determine suitable sampling densities $K_L$ and $K_\theta$ and energy fraction tolerance $\delta_{\text{EF}}$ for accurate ROM-MD simulations, 
we consider a force difference of $5 \times 10^{-4}$ Hartree/Bohr as a negligible error tolerance. 
We start by evaluating the performance of our reduced basis on a single, unseen molecular configuration 
defined by $s_1 = 1.02$, $s_2 = 0.98$, and $s_\theta = 2.0$. This testing configuration lies within the parametric domain but was not included in the training set.
Table~\ref{tab:distance_to_training} shows the distance between this testing configuration and 
the closest atomic system in the training set in the parametric space for different sampling densities. 
The distances for $s_1$ and $s_2$ are calculated as the absolute difference between the testing value and the nearest sampled $s$ value, multiplied by the scaling factor 1.83 Bohr. 
The distance for $s_\theta$ is the absolute difference in degrees. 

\begin{table}[h!]
    \centering
    \begin{tabular}{|c||cccc|}
        \hline
        \multirow{2}{*}{Parameter} & \multicolumn{4}{c|}{$(K_L, K_\theta)$} \\
        & $(2,2)$ & $(2,5)$ & $(5,2)$ & $(5,5)$ \\ 
        \hline
        \hline
        $L_1$ & 0.0366 & 0.0366 & 0.0183 & 0.0183 \\
        $L_2$ & 0.0366 & 0.0366 & 0.0183 & 0.0183 \\
        $\theta$ & 2.0 & 1.0 & 2.0 & 1.0 \\
        \hline
    \end{tabular}
    \caption{Distance to the closest training atomic system in parametric space for the testing configuration ($s_1 = 1.02, s_2 = 0.98, s_\theta = 2.0$).}
    \label{tab:distance_to_training}
\end{table}

Tables~\ref{tab:force_error_H1} and~\ref{tab:force_error_H2} present the magnitude of the force difference on hydrogen atom \ce{H1} and \ce{H2}, respectively, 
between the full first-principles calculation and the data-driven ROM for the unseen testing configuration. 
These force differences are evaluated for different sampling densities and energy fraction tolerance $\delta_{\text{EF}}$.
The results indicate that increasing the sampling density and decreasing the energy fraction tolerance $\delta_{\text{EF}}$ 
generally leads to a smaller force difference, as expected from the increased accuracy of the reduced basis. 
Notably, for the lowest sampling frequency $K_L = K_\theta = 2$ and energy fraction residual $\delta_{\text{EF}} = 10^{-4}$, 
the force differences on both hydrogen atoms are consistently below our chosen tolerance of $5 \times 10^{-4}$ Hartree/Bohr. 
For the remaining studies, we will use the subspace of dimension $r = 34$ generated from this combination of sampling frequency and energy fraction residual. 

\begin{table}[h!]
    \centering
    \begin{tabular}{|c||cccc|}
        \hline
        \multirow{2}{*}{$\delta_{\text{EF}}$} & \multicolumn{4}{c|}{$(K_L, K_\theta)$} \\
        & $(2,2)$ & $(2,5)$ & $(5,2)$ & $(5,5)$ \\ 
        \hline
        \hline
        $10^{-2}$ & $7.3 \times 10^{-4}$ & $2.6 \times 10^{-2}$ & $2.1 \times 10^{-2}$ & $2.1 \times 10^{-2}$ \\
        $10^{-3}$ & $2.0 \times 10^{-4}$ & $1.8 \times 10^{-4}$ & $1.5 \times 10^{-4}$ & $1.4 \times 10^{-4}$ \\
        $10^{-4}$ & $2.4 \times 10^{-5}$ & $3.1 \times 10^{-5}$ & $2.4 \times 10^{-5}$ & $2.3 \times 10^{-5}$ \\
        \hline
    \end{tabular}
    \caption{Magnitude of the difference in force on hydrogen atom \ce{H1} (in Hartree/Bohr) between FPMD and ROM-MD for the testing configuration.}
    \label{tab:force_error_H1}
\end{table}

\begin{table}[h!]
    \centering
    \begin{tabular}{|c||cccc|}
        \hline
        \multirow{2}{*}{$\delta_{\text{EF}}$} & \multicolumn{4}{c|}{$(K_L, K_\theta)$} \\
        & $(2,2)$ & $(2,5)$ & $(5,2)$ & $(5,5)$ \\ 
        \hline
        \hline
        $10^{-2}$ & $4.8 \times 10^{-4}$ & $9.6 \times 10^{-3}$ & $2.1 \times 10^{-2}$ & $2.1 \times 10^{-2}$ \\
        $10^{-3}$ & $1.0 \times 10^{-4}$ & $7.2 \times 10^{-5}$ & $5.2 \times 10^{-5}$ & $7.0 \times 10^{-5}$ \\
        $10^{-4}$ & $4.3 \times 10^{-5}$ & $4.6 \times 10^{-5}$ & $3.9 \times 10^{-5}$ & $3.8 \times 10^{-5}$ \\
        \hline
    \end{tabular}
    \caption{Magnitude of the difference in force on hydrogen atom \ce{H2} (in Hartree/Bohr) between FPMD and ROM-MD for the testing configuration.}
    \label{tab:force_error_H2}
\end{table}

\begin{figure}[h!]
    \centering
    \includegraphics[width=0.49\textwidth]{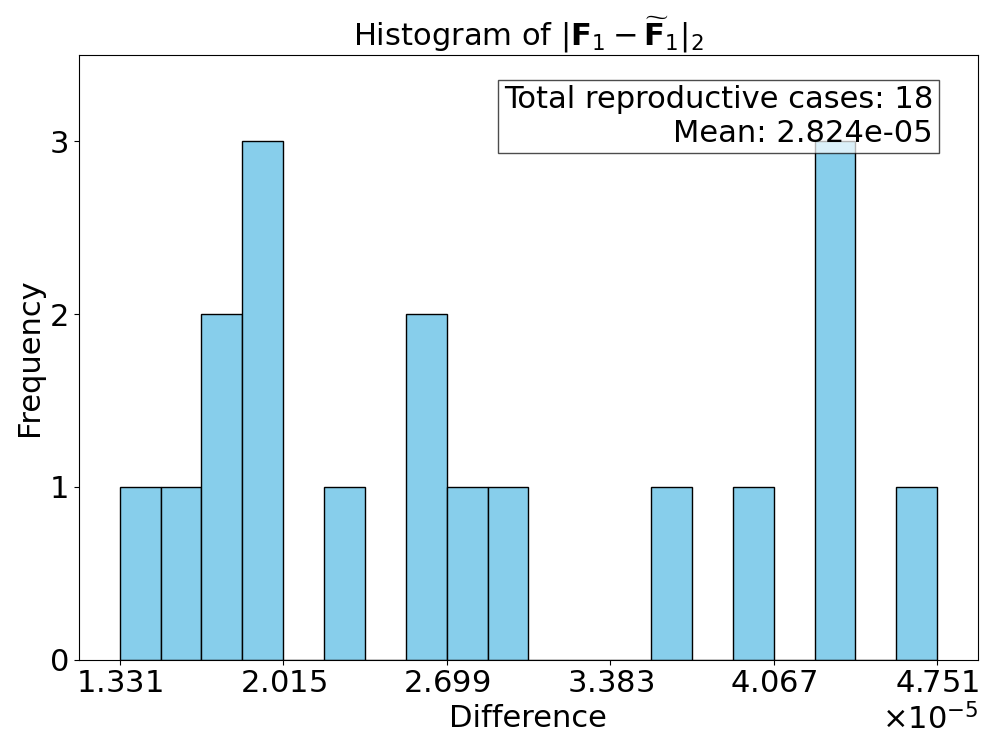}
    \includegraphics[width=0.49\textwidth]{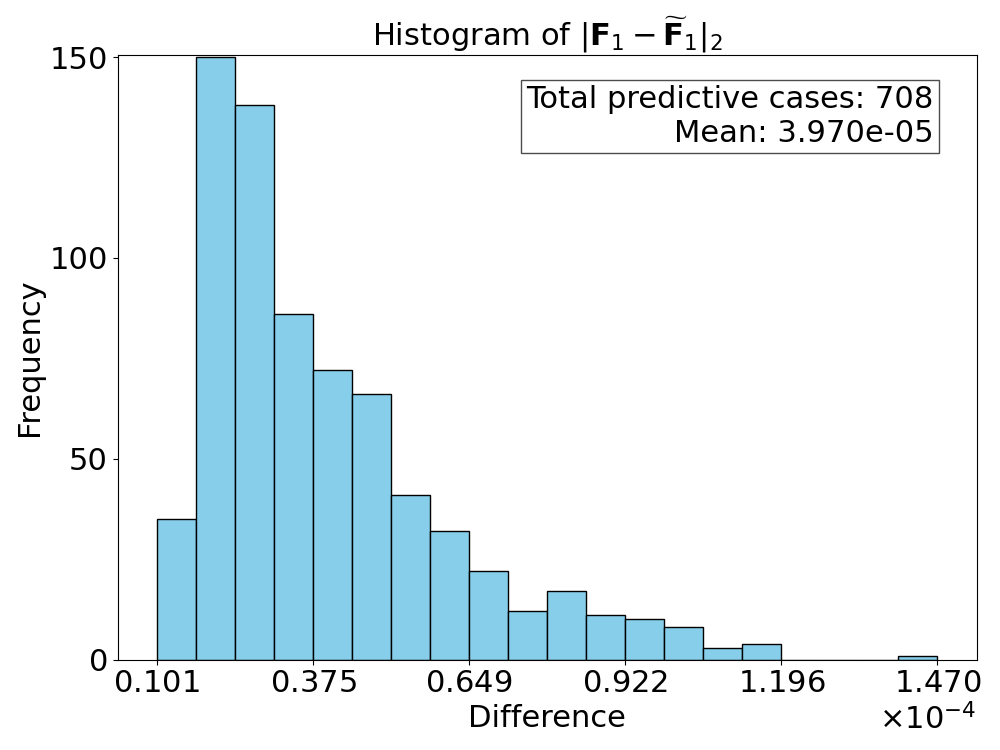} \\
    \includegraphics[width=0.49\textwidth]{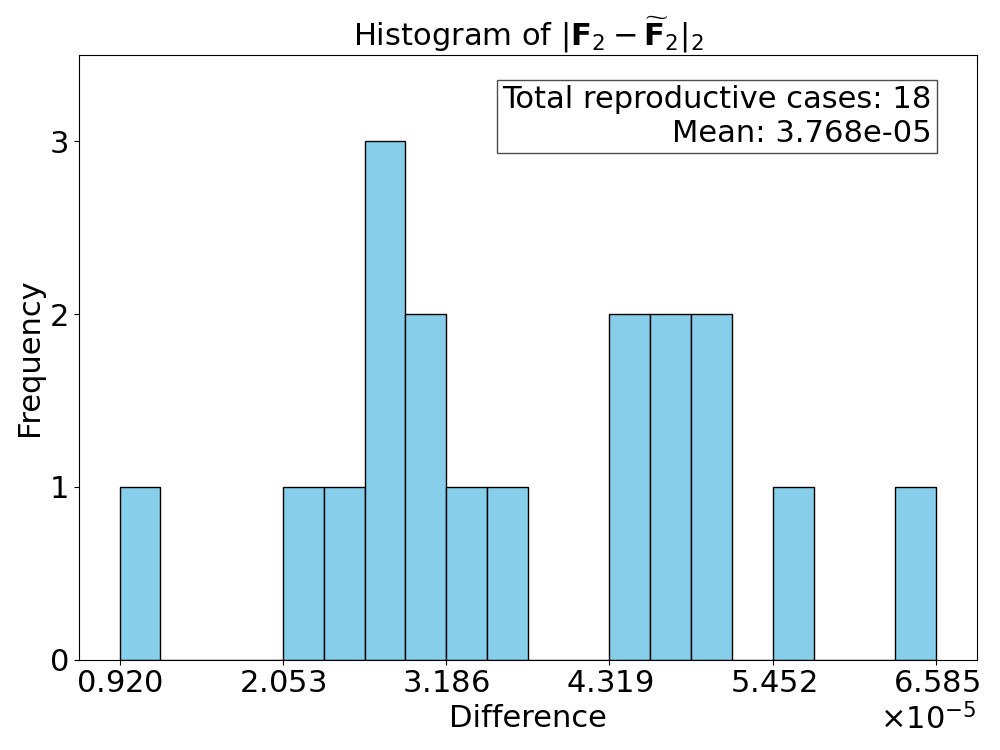}
    \includegraphics[width=0.49\textwidth]{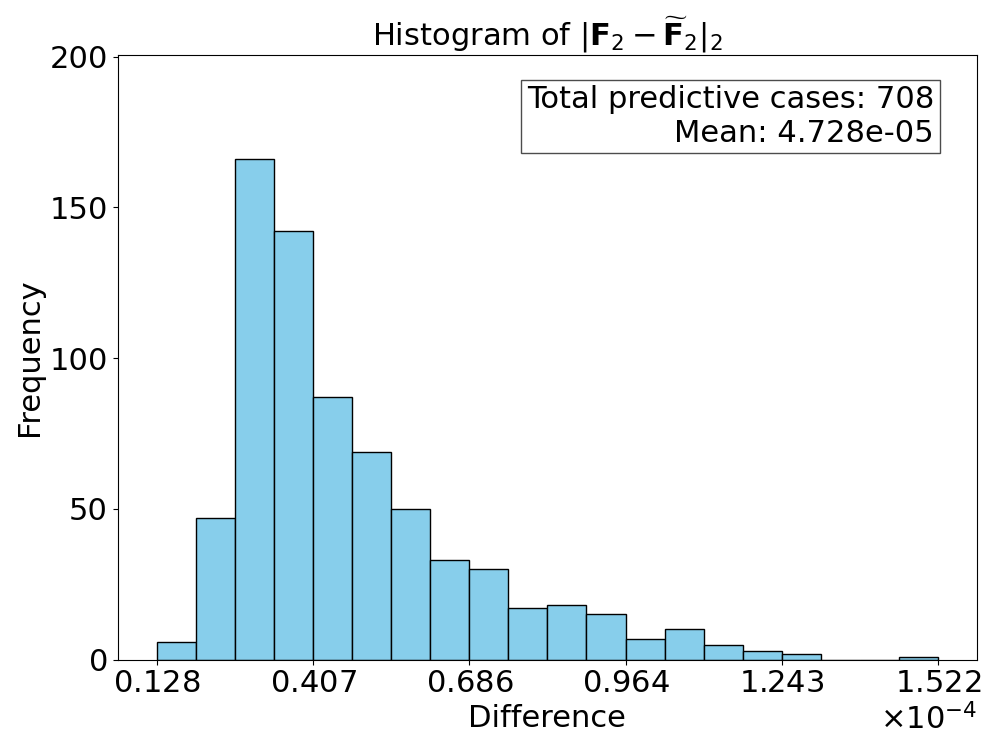}
    \caption{Histogram of the magnitude of the difference in force (in Hartree/Bohr) between FPMD and ROM-MD for hydrogen atoms \ce{H1} (top row) and \ce{H2} (bottom row). The left column shows results for the 18 reproductive cases, and the right column shows results for the 708 predictive unseen configurations. The mean force difference is reported in each panel.}
    \label{fig:force_histograms}
\end{figure}

Following the single-point validation, we further assessed the performance of our ROM-MD approach by evaluating the force differences 
on the hydrogen atoms (\ce{H1} and \ce{H2}) over a larger set of $K_{\text{test}}=726$ atomic configurations. 
These configurations were generated by uniformly sampling the parametric space with a finer grid with $K_L = K_\theta = 10$, 
ensuring a diverse set of molecular geometries within our defined domain. 
Among these tested configurations, 18 are reproductive cases which coincided exactly with configurations present in our training set, 
while the remaining 708 were unseen predictive configurations. 
The histograms in Figure~\ref{fig:force_histograms} reveal the distribution of the magnitude of force differences for both hydrogen atoms across the reproductive and predictive test sets. 
For the 18 reproductive configurations (left column), the mean force difference is $2.9520 \times 10^{-5}$ Hartree/Bohr for \ce{H1} and $3.3766 \times 10^{-5}$ Hartree/Bohr for \ce{H2}, 
and the maximum difference is $6.767 \times 10^{-5}$ Hartree/Bohr for \ce{H1} and $6.612 \times 10^{-5}$ Hartree/Bohr for \ce{H2}. 
For the 708 unseen predictive configurations (right column), the distribution of force differences exhibits a right-skewed shape with some outliers. 
The mean force difference for \ce{H1} in the predictive set is $3.950 \times 10^{-5}$ Hartree/Bohr, and for \ce{H2}, the mean is $6.734 \times 10^{-5}$ Hartree/Bohr. 
Importantly, despite the skewed distribution and the presence of outliers, the maximum force difference $1.670 \times 10^{-4}$ Hartree/Bohr for \ce{H1} and $1.521 \times 10^{-4}$ Hartree/Bohr for \ce{H2} 
in the predictive set remains below our defined negligible force difference threshold of $5 \times 10^{-4}$ Hartree/Bohr. 
This suggests that our ROM-MD approach exhibits good predictive capability for unseen configurations within the sampled parametric domain, maintaining force accuracy within an acceptable tolerance.

\subsection{Accuracy in atomic trajectories}

To empirically validate our ROM-MD framework proposed in Section~\ref{sec:rom}, 
we present an example simulation of a water molecule. 
The initial geometry of the water molecule is provided in Table~\ref{tab:initial_geometry}. 

\begin{table}[h!]
    \centering
    \begin{tabular}{|c||c c c|}
        \hline
        Atom & $x$ & $y$ & $z$ \\
        \hline
        \hline
        \ce{O} & 0.00 & 0.00 & 0.00 \\
        \ce{H1} & -0.45 & 1.57 & -1.07 \\
        \ce{H2} & -0.45 & -1.48 & -0.97 \\
        \hline
    \end{tabular}
    \caption{Initial geometry of the water molecule in Bohr.}
    \label{tab:initial_geometry}
\end{table}

\begin{figure}[h!]
    \centering
    \includegraphics[width=0.6\linewidth]{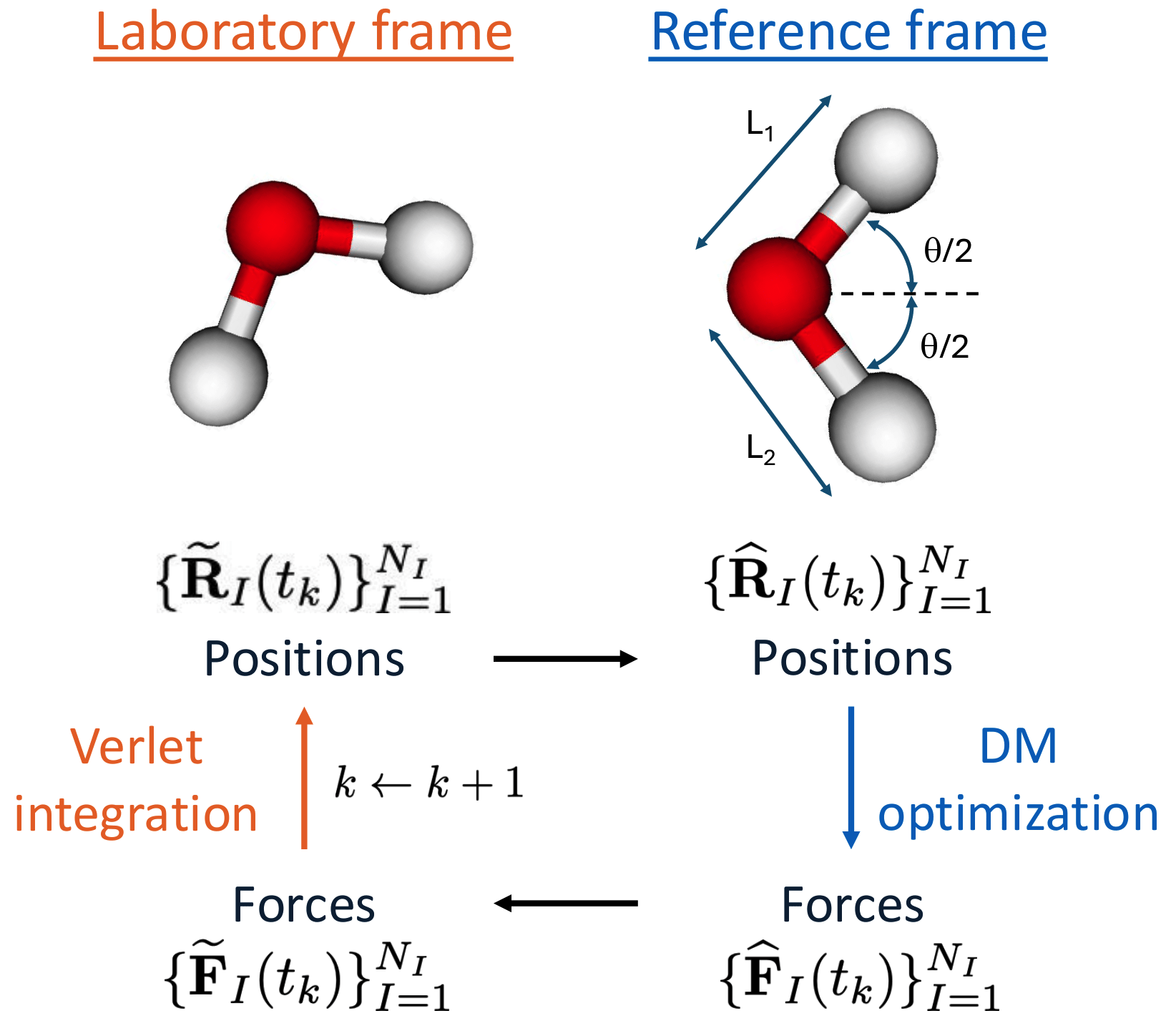}
    \caption{Schematic diagram for the procedure of reduced Born-Oppenheimer molecular dynamics in Algorithm~\ref{alg:rom-born-oppenheimer} 
in the example of water molecule.}
    \label{fig:schematic}
\end{figure}

The molecular dynamics simulation was run for $K_{\text{MD}} = 500$ time steps with a time step of 40.0 atomic units (approximately 0.967 femtoseconds). 
No thermostat was used for this simulation, allowing us to check for energy conservation. Given the small problem size, large temperature fluctuations are expected.
We remark that the initial atomic system, with the given coordinates, does not lie exactly within our previously defined parametric domain. 
This choice was made to test the extrapolative capabilities of our ROM-MD approach 
beyond the configurations explicitly included in the training data. 
By slightly perturbing the initial positions away from the equilibrium configuration, we initiate a dynamic trajectory where the temperature will naturally evolve.
The temperature remains within a range of 0 to 400 Kelvin during the simulation 
avoiding excessive extrapolation beyond the training data's energy landscape.
Figure~\ref{fig:schematic} illustrates the procedure of reduced Born-Oppenheimer molecular dynamics in Algorithm~\ref{alg:rom-born-oppenheimer} 
in the example of water molecule. 

The performance of the ROM-MD method is significantly affected by the choice of the convergence tolerance $\delta_{\text{DM}}$ in Algorithm~\ref{alg:dm-solver}. Therefore, a systematic study of the relative speed-up as a function of these parameters is necessary. 
With a mixing parameter $\beta = 0.8$, 
the relative speed-up of the ROM-MD simulation compared to the FPMD simulation is presented in Table~\ref{tab:speed-up}. The speed-up is defined as the ratio of the total wall-clock time of the FPMD simulation to that of the ROM-MD simulation. From the table, we observe that for less stringent tolerances, such as $\delta_{\text{DM}} = 10^{-4}$, the ROM-MD method achieves a significant speed-up of over 4 times. However, as the tolerance is tightened to $\delta_{\text{DM}} = 10^{-8}$, the speed-up decreases. 
This is because the evaluation of the nonlinear terms is a significant overhead, and this overhead increases with a tighter tolerance. 
The results presented in the rest of the paper will use the parameters of $\beta = 0.8$ and $\delta_{\text{DM}} = 10^{-6}$, 
which represent a balance between computational efficiency and accuracy.

\begin{table}[h!]
\centering
\begin{tabular}{|c||c c c|} 
\hline
$\delta_{\text{DM}}$ & $10^{-4}$ & $10^{-6}$ & $10^{-8}$ \\ 
\hline
\hline
Factor & 4.17 & 2.05 & 1.39 \\
\hline
\end{tabular}
\caption{Relative speed-up of ROM-MD with different choices of 
convergence tolerance $\delta_{\text{DM}}$ in Algorithm~\ref{alg:dm-solver}.}
\label{tab:speed-up}
\end{table}

\begin{figure}[h!]
    \centering
    \includegraphics[width=0.49\textwidth]{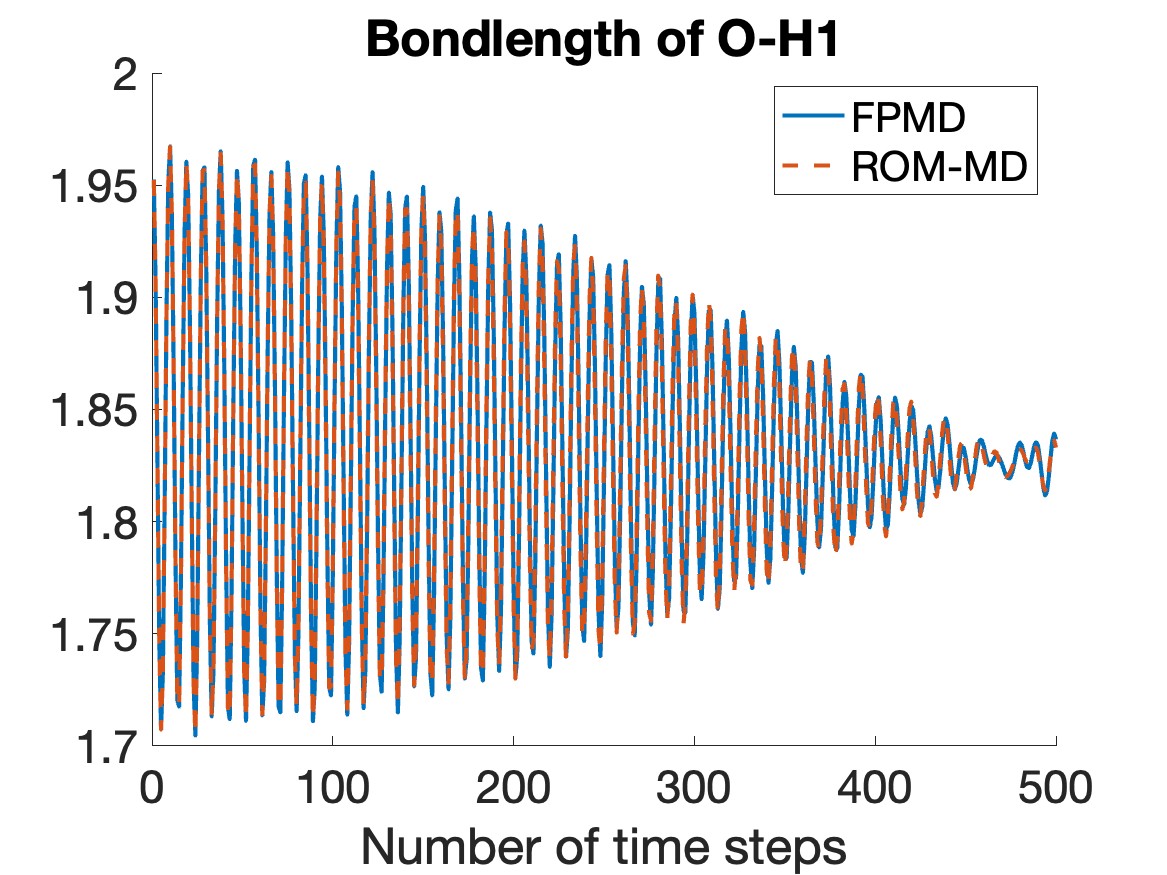}
    \includegraphics[width=0.49\textwidth]{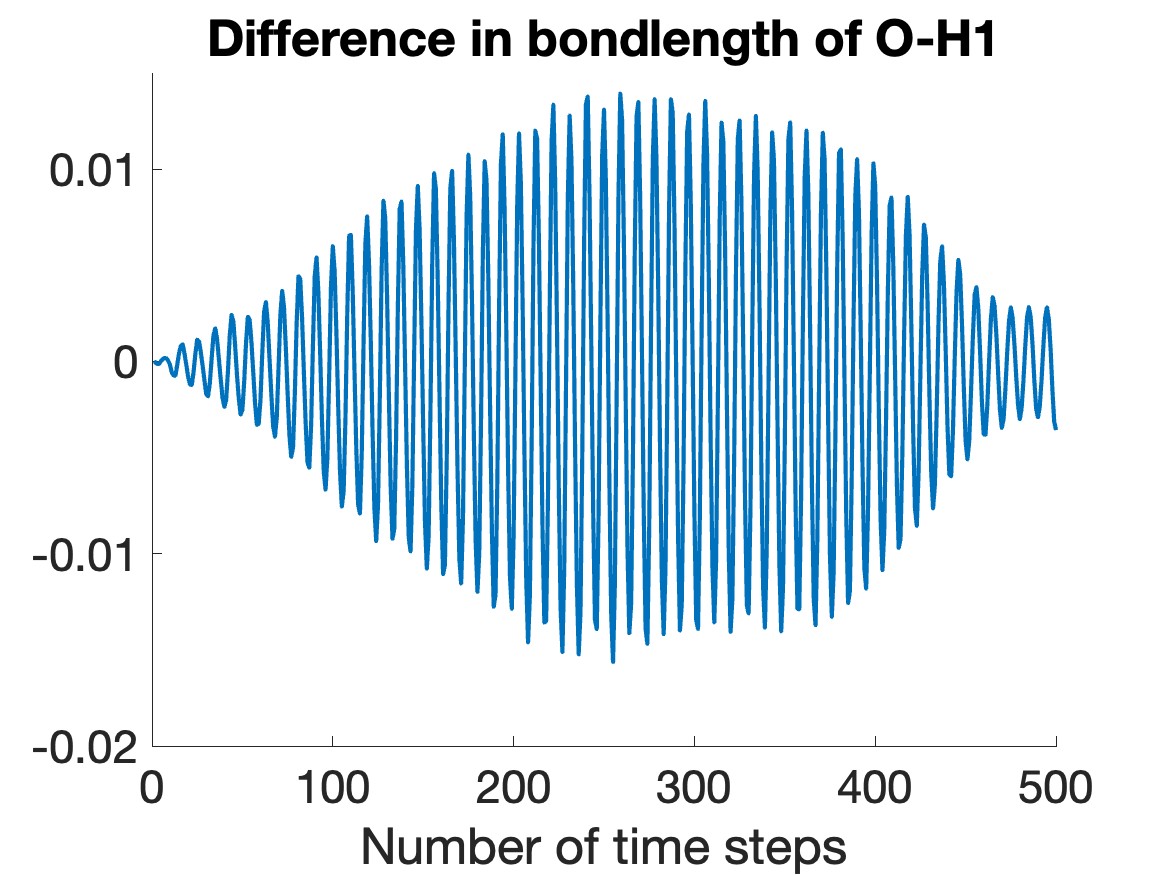}\\
    \includegraphics[width=0.49\textwidth]{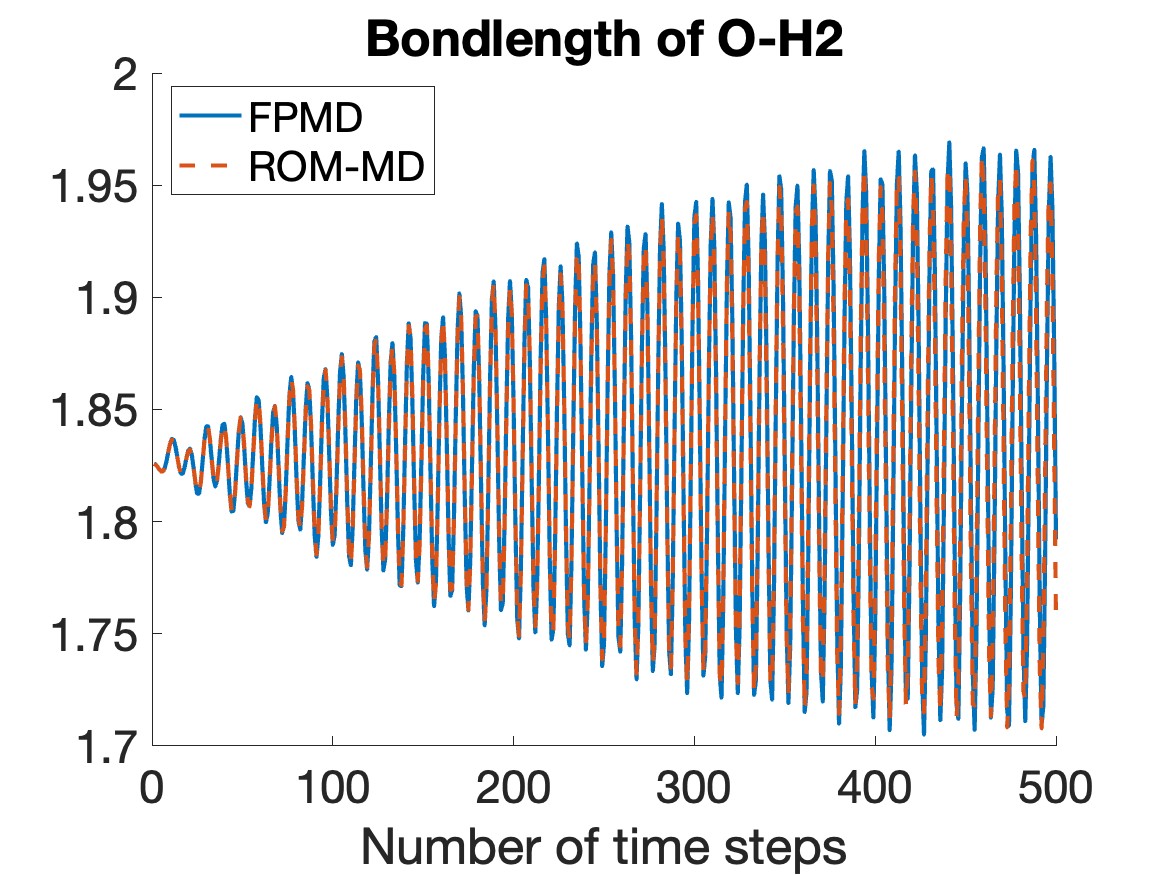}
    \includegraphics[width=0.49\textwidth]{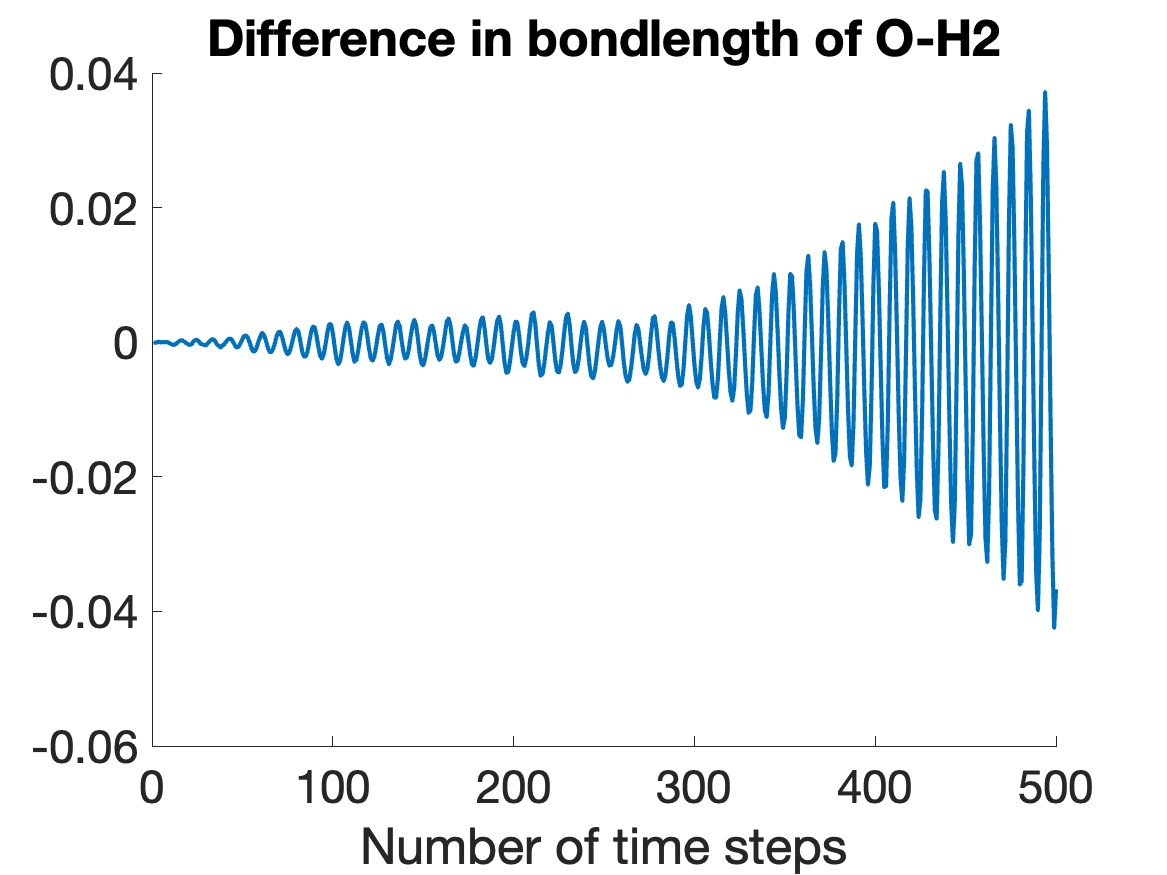}\\
    \includegraphics[width=0.49\textwidth]{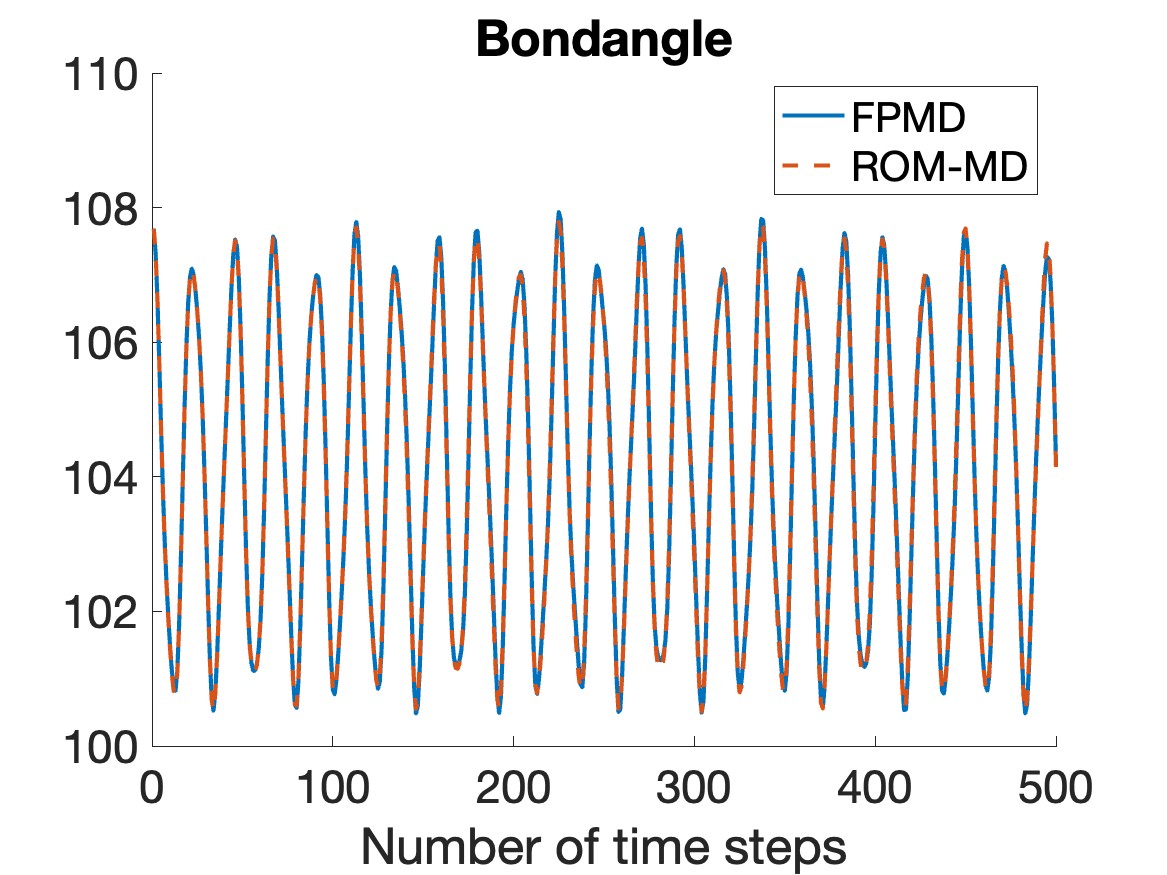}
    \includegraphics[width=0.49\textwidth]{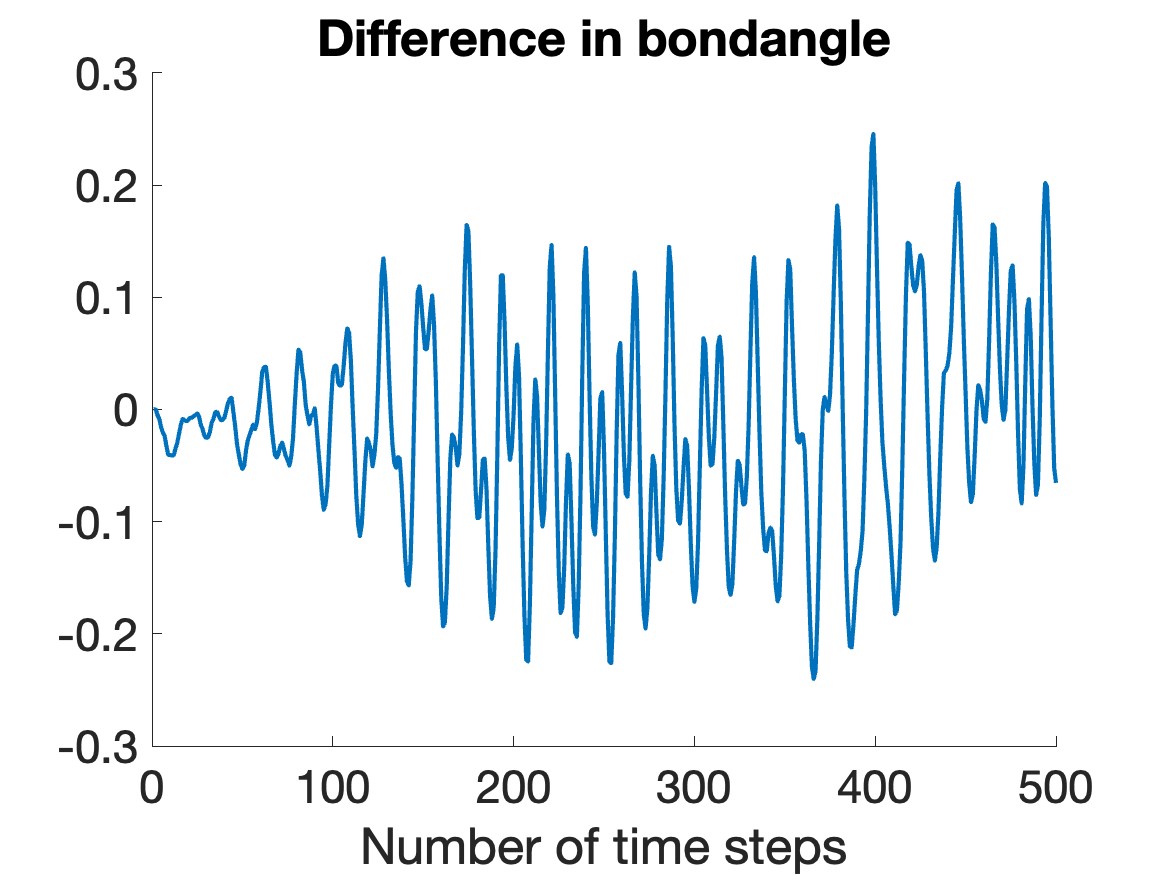}\\ 
    \caption{Evolution of bond lengths $L_1$ and $L_2$ and bond angle $\theta$ during the 500 time steps of molecular dynamics. 
    Left column: Comparison between FPMD (blue) and ROM-MD (red). Right column: Difference between FPMD and ROM-MD results.}
    \label{fig:bond_properties_1e-6}
\end{figure}

Figure~\ref{fig:bond_properties_1e-6} shows the evolution of key geometric variables of the water molecule during the 500 time steps of the molecular dynamics simulation, 
comparing the FPMD simulation in Section~\ref{sec:fom} and the ROM-MD simulation in Section~\ref{sec:rom}. 
The left column shows the length $L_1$ of the bond \ce{O-H1} (left panel), 
$L_2$ of the bond \ce{O-H2} (middle panel), 
and the bond angle $\theta$ formed by the two \ce{O-H} bonds (right panel).
In each of these top panels, the results obtained from the FPMD simulation and the ROM-MD simulation are shown in blue and red respectively. 
The plots reveal that the ROM-MD simulation resemble the FPMD simulation extremely well and capture the oscillatory and periodic nature of the hydrogen atoms' motion relative to the fixed oxygen. 
The red curves from the ROM-MD simulation closely follow the blue curves from the FPMD simulation for all three geometric properties. 
The right column, illustrating the differences between the two methods, shows that the deviations remain relatively small throughout the simulation, 
indicating a good agreement between the FPMD approach and ROM-MD approach in predicting the molecular geometry over time.

\begin{figure}[h!]
    \centering
    \includegraphics[width=0.49\textwidth]{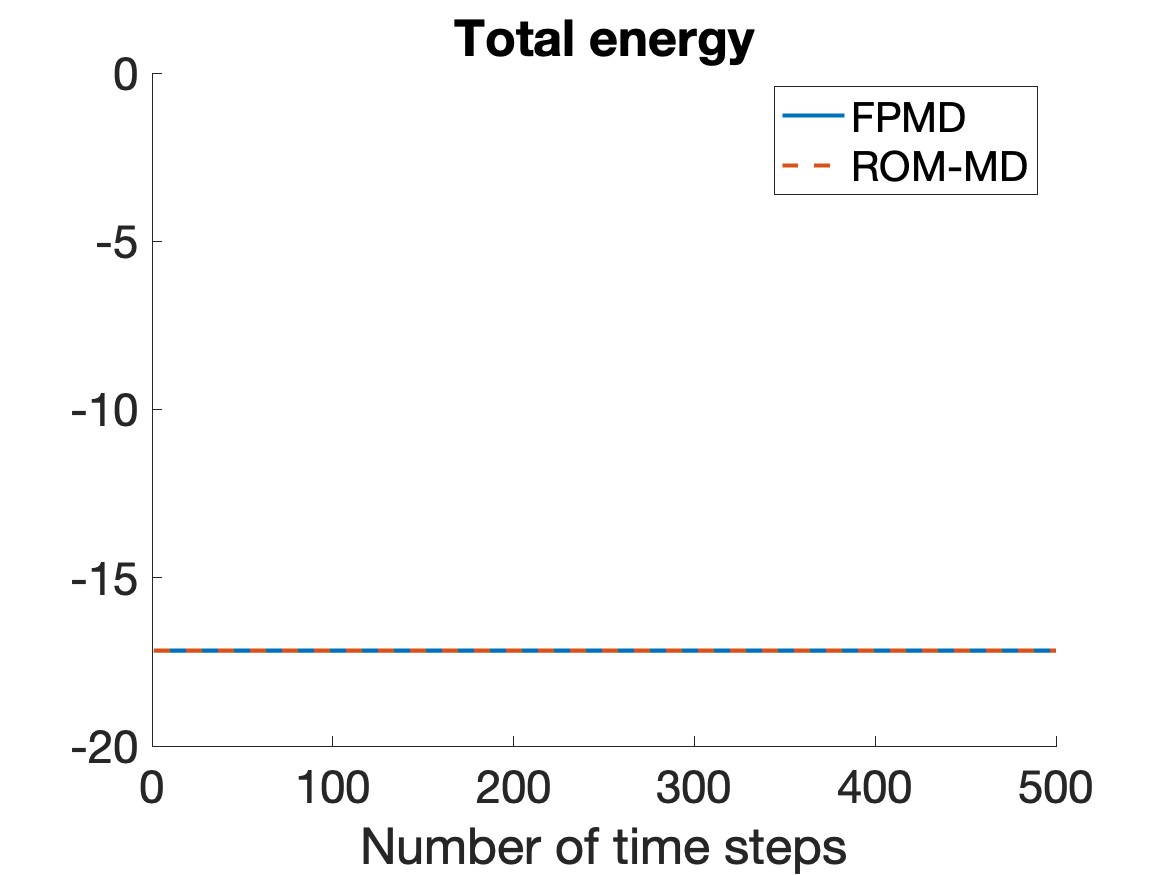}
    \includegraphics[width=0.49\textwidth]{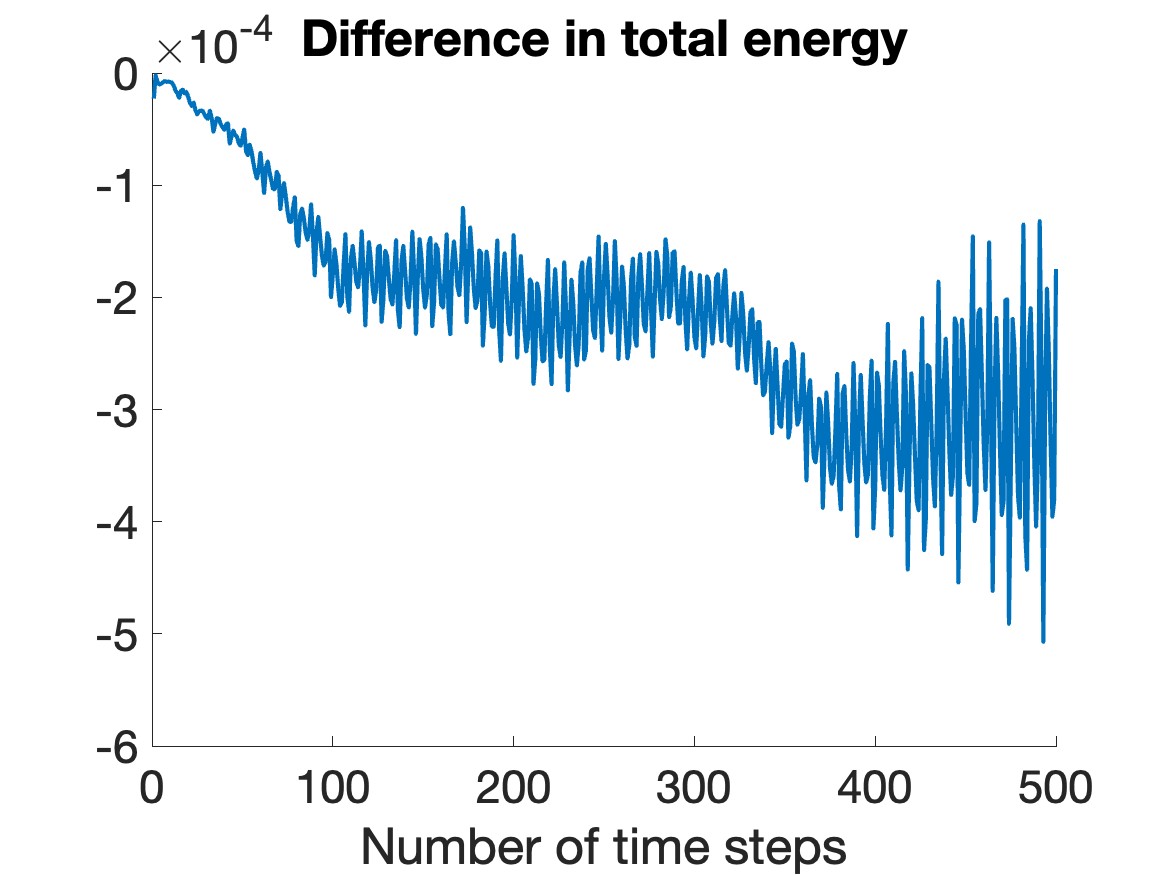}\\ 
    \caption{Evolution of total energy during the 500 time steps of molecular dynamics. 
    Left column: Comparison between FPMD (blue) and ROM-MD (red). Right column: Difference between FPMD and ROM-MD results.}
    \label{fig:energy_conservation_1e-6}
\end{figure}

Figure~\ref{fig:energy_conservation_1e-6} illustrates the evolution of total energy of the system during the 500 time steps of the molecular dynamics simulation, 
comparing the FPMD simulation in Section~\ref{sec:fom} and the ROM-MD simulation in Section~\ref{sec:rom}. 
Both the FPMD and ROM-MD methods successfully conserve the total energy of the system. The total energy plots for both methods fluctuate within a very narrow range 
centered around -17.165 Hartree. This small fluctuation demonstrates that the simulations are stable and accurately maintain energy conservation throughout the 500 time steps. The trajectory of the difference further confirms the high degree of similarity between the two methods by plotting the difference in their total energies. This difference remains extremely small, on the order of $10^{-4}$, which indicates that ROM-MD simulation has a comparable total energy to the full-scale FPMD simulation throughout 500 time steps.

%% file: conclusion.tex
\section{Conclusion}

In this work, we have presented a data-driven framework for ground state electronic calculations based on learning a reduced representation of the electronic wavefunctions. 
Through a systematic sampling of the configuration space of atomic coordinates and the construction of a low-dimensional basis for the electronic wavefunctions using snapshot singular value decomposition, 
we have shown the feasibility of approximating interatomic forces with excellent accuracy for a pinned water molecule. This reduction successfully captured the essential physics of the electronic wave function space with a significantly smaller basis.

When applied to ROM-MD, the atomic forces approximated using this reduced representation closely reproduced the bond lengths and bond angle obtained from full first-principles calculations. The accuracy of these predicted molecular properties serves as a powerful validation of the quality of the learned reduced electronic space.
%While our initial implementation had not produced significant computational speed-ups, primarily due to the direct evaluation of nonlinear terms, 
%the accuracy of the predicted molecular properties underscores the potential of this approach. 
To unlock the full potential of the low dimensional representation for achieving computational speed-ups, 
an interesting direction is the use of data-driven surrogates 
to update effective Hamiltonian, which is crucial for achieving substantial acceleration and enabling the simulation of larger and more complex systems.
An example is hyper-reduction techniques \cite{chaturantabut2010nonlinear,carlberg2013gnat,drmac2016new,drmac2018discrete,lauzon2024s} 
to efficiently compute these nonlinear components within the reduced basis of electronic density.

Another direction for further study, which we begin exploring and report preliminary results in Appendix~\ref{appendix}, is the use of nonlinear compression techniques to obtain even lower-dimensional, highly-accurate representations of wavefunction solutions.
This shift investigates a more general, nonlinear reduction of the electronic manifold, which is the key goal of advanced ROM for electronic structure.
This study provides a foundational step towards developing nonlinear ROM for FPMD simulations.

\section*{Acknowledgement}
This work was supported by Laboratory Directed Research and Development (LDRD) Program by the U.S. Department of Energy (24-ERD-035). 
Lawrence Livermore National Laboratory is operated by Lawrence Livermore National Security, LLC, for the U.S. Department of Energy, National Nuclear Security Administration under Contract DE-AC52-07NA27344. IM release number: LLNL-JRNL-2013760.
This manuscript has been co-authored by UT-Battelle, LLC, under contract DE-AC05–00OR22725 with the US Department of Energy (DOE).

%% file: appendix.tex
\appendix
\section{Nonlinear compression} \label{appendix}

To explore the lowest-dimensional structures that may accurately represent electronic wavefunctions, we turn to nonlinear compression.
Here, we test an autoencoder (AE) approach for the same pinned water molecule example considered in Section \ref{sec:numerical}.
Taking advantage of the fact that linear compression already performs quite well for this problem, we employ a learnable weighted hybrid AE that combines a (linear) POD-based AE with a (nonlinear) convolutional AE \cite{somasekharan2025beyond}.
Each AE $\mathbf{A}$ consists of an encoder $\mathbf{E}$ that maps an input vector $\mathbf{x} \in \mathbb{R}^{M}$ to a lower-dimensional latent vector $\mathbf{z} \in \mathbb{R}^{r}$ and a decoder $\mathbf{D}$ that maps the latent vector to $\tilde{\mathbf{x}} \in \mathbb{R}^{M}$, such that $\tilde{\mathbf{x}} = \mathbf{A} \left( \mathbf{x} \right) = \mathbf{D} \left( \mathbf{E} \left( \mathbf{x} \right) \right)$ is the reconstruction of the compressed input.
In the POD case, $\mathbf{E}_{\mathrm{POD}}\left( \mathbf{x} \right) = \mathbf{Q}^\top \mathbf{x}$ and $\mathbf{D}_{\mathrm{POD}}\left( \mathbf{z} \right) = \mathbf{Q} \mathbf{z}$, where $\mathbf{Q} \in \text{St}(M, r)$ is the ROM basis matrix of rank $r$ defined in Section \ref{sec:offline}. The hybrid AE combines these with a corresponding convolutional neural network-based encoder $\mathbf{E}_{\mathrm{CNN}}$ and decoder $\mathbf{D}_{\mathrm{CNN}}$:

\begin{equation}
    \mathbf{E}_\mathrm{hybrid} \left( \mathbf{x} \right)
    = \left( \mathbf{1} - \mathbf{a} \right) \odot \mathbf{E}_{\mathrm{POD}}\left( \mathbf{x} \right)
    + \mathbf{a} \odot \mathbf{E}_{\mathrm{CNN}}\left( \mathbf{x} \right)
    ,
\end{equation}

\begin{equation}
    \mathbf{D}_\mathrm{hybrid} \left( \mathbf{z} \right) = \left( \mathbf{1} - \mathbf{b} \right) \odot \mathbf{D}_{\mathrm{POD}}\left( \mathbf{z} \right)
    + \mathbf{b} \odot \mathbf{D}_{\mathrm{CNN}}\left( \mathbf{z} \right)
    .
\end{equation}
$\mathbf{a} \in \mathbb{R}^{r}$ and $\mathbf{b} \in \mathbb{R}^{M}$ here are learnable weight vectors, while $\odot$ denotes element-wise multiplication.

In this case, we use a reduce dimension of $r=18$ in the bottleneck layer of the AE.
The convolutional encoder $\mathbf{E}_{\mathrm{CNN}}$ consists of four 3D convolution layers with decreasing spatial dimension (32 to 16 to 8 to 4) and increasing channels (32 to 64 to 128 to 256) in their outputs, followed by a single fully-connected layer that reduces the total dimension from 16384 to 18.
The convolutional decoder $\mathbf{D}_{\mathrm{CNN}}$ mirrors the encoder, with one fully-connected layer followed by four 3D transposed convolution layers.
ReLU activation functions are applied to the output of each layer except for the final layer of the decoder, for which we use the hyperbolic tangent function.
As in \cite{somasekharan2025beyond}, we choose to initialize the weights $\mathbf{a}$ and $\mathbf{b}$ to zero, such that $\mathbf{A}_\mathrm{hybrid}$ is initially equivalent to $\mathbf{A}_\mathrm{POD}$.
For training, we use the same sampling frequency of molecular configurations ($K_L = K_\theta = 2$) and divide the 72 wavefunctions into four batches of 18.
The hybrid AE is implemented and trained using PyTorch \cite{paszke2019pytorch}.
We use the Adam optimizer \cite{kingma2017adammethodstochasticoptimization} with a fixed learning rate of $5 \times 10^{-4}$ and $\epsilon = 1 \times 10^{-12}$ and train for 200000 epochs.

Table \ref{tab:force_error_ae} lists the resulting relative error $\frac{\vert \mathbf{Y} - \widetilde{\mathbf{Y}} \vert_\text{F}}{\vert \mathbf{Y} \vert_\text{F}}$ in the snapshot matrix $\mathbf{Y}$ for the POD-based and hybrid AEs (where $\widetilde{\mathbf{Y}}$ is the AE-reconstructed snapshot matrix).
Note that due to the properties of the Frobenius norm, this error is equivalent to $\delta_{\text{EF}}$ in the POD case.
The hybrid approach reduces the error by a factor of more than 12 relative to POD alone, such that it is below $10^{-4}$.
Reaching this target energy fraction residual requires $r=34$ in the linear ROM case (Table \ref{tab:rom_dimension}) but is achieved here with $r=18$.

\begin{figure}[h!]
    \centering
    \includegraphics[width=0.49\textwidth]{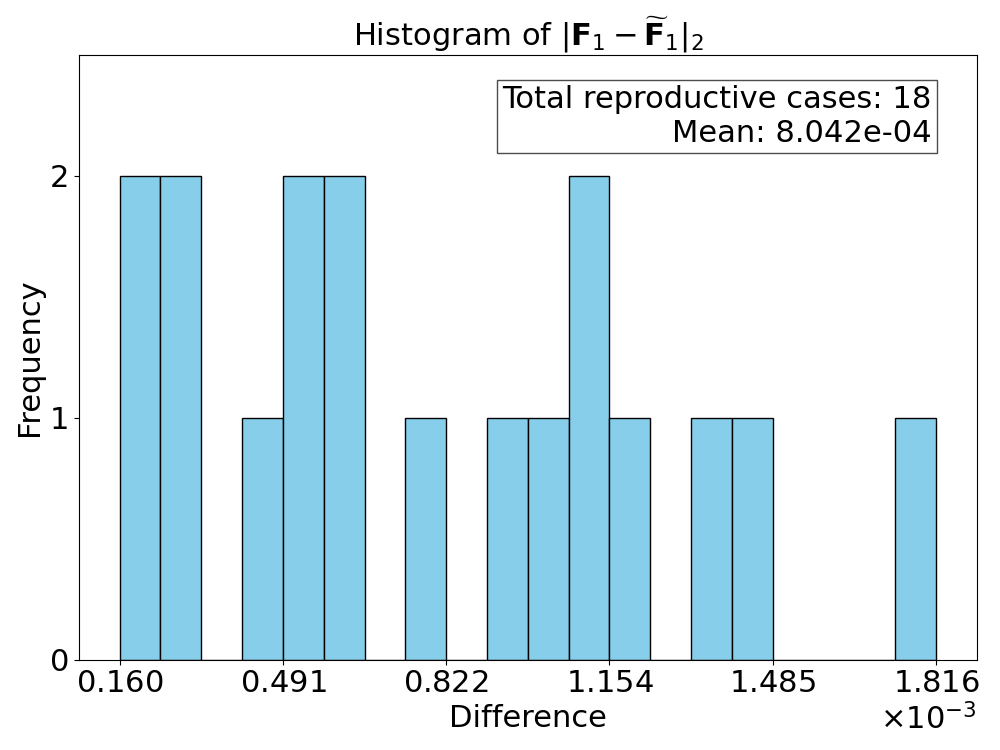}
    \includegraphics[width=0.49\textwidth]{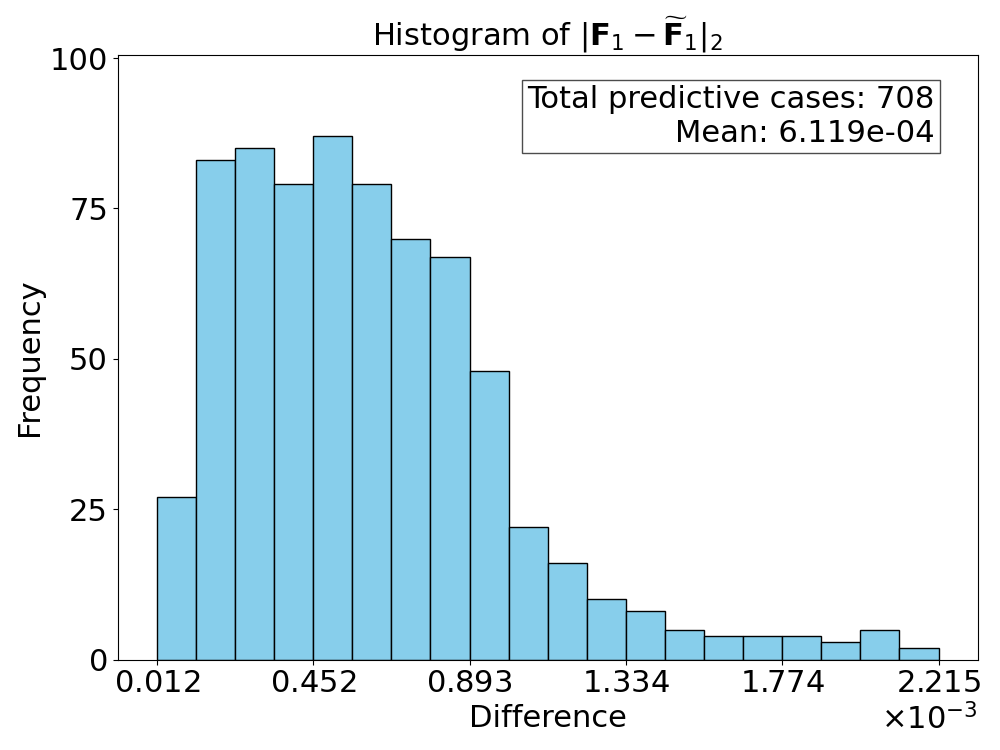} \\
    \includegraphics[width=0.49\textwidth]{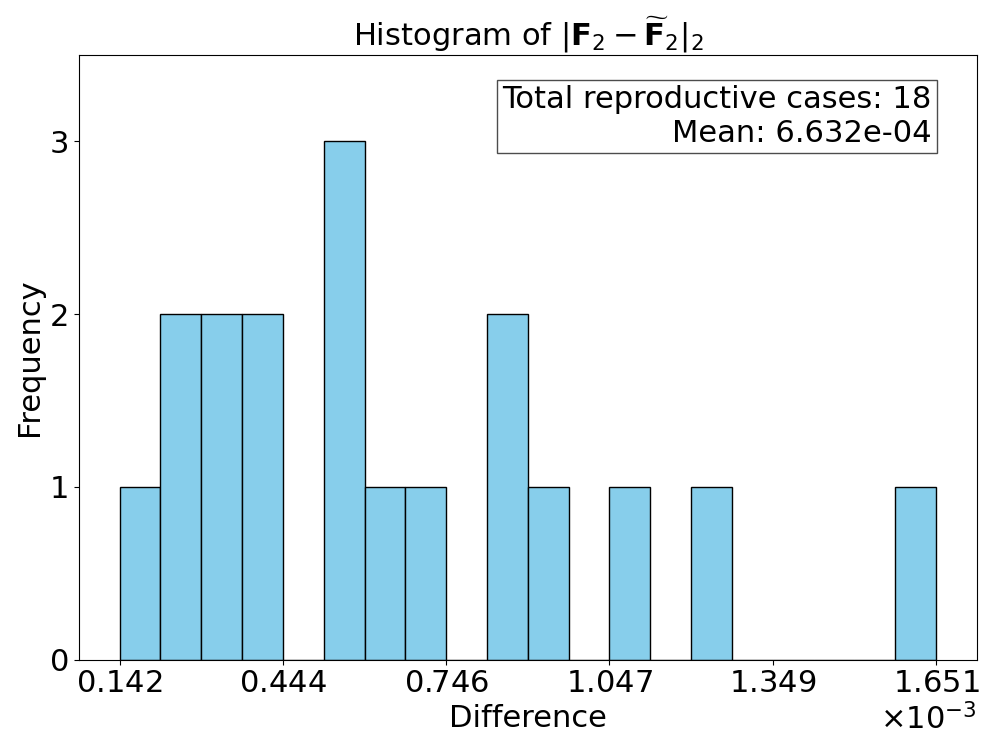}
    \includegraphics[width=0.49\textwidth]{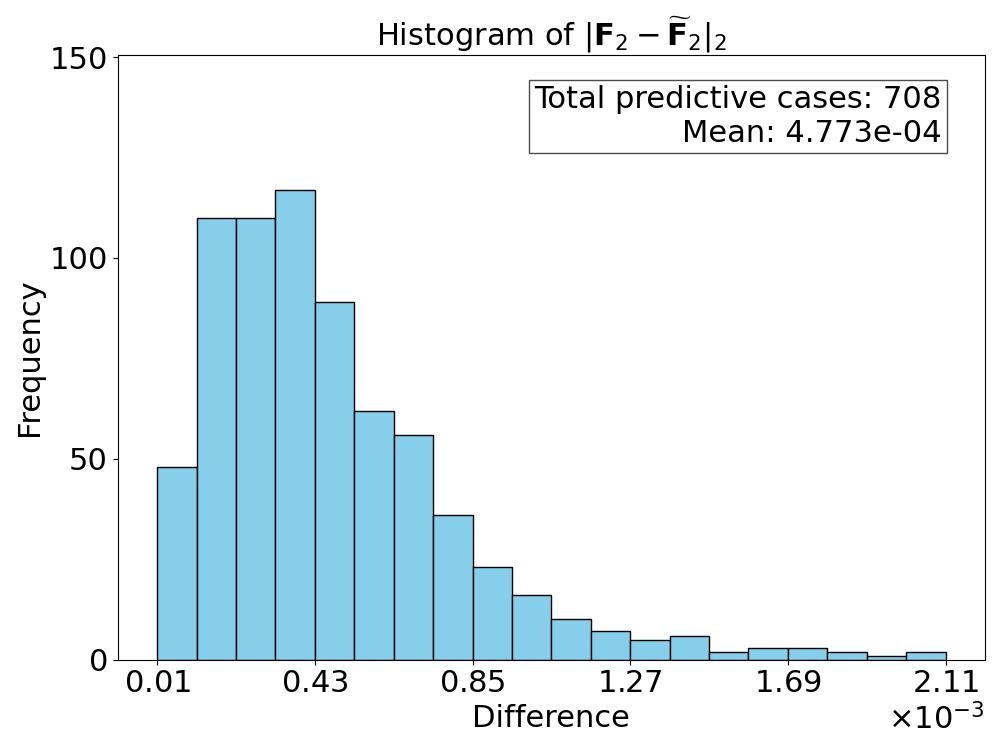}
    \caption{Histogram of the magnitude of the difference in force (in Hartree/Bohr) between FPMD and POD-based AE compression for hydrogen atoms \ce{H1} (top row) and \ce{H2} (bottom row). The left column shows results for the 18 reproductive cases, and the right column shows results for the 708 predictive unseen configurations. The mean force difference is reported in each panel.}
    \label{fig:force_histograms_pod_18}
\end{figure}

\begin{figure}[h!]
    \centering
    \includegraphics[width=0.49\textwidth]{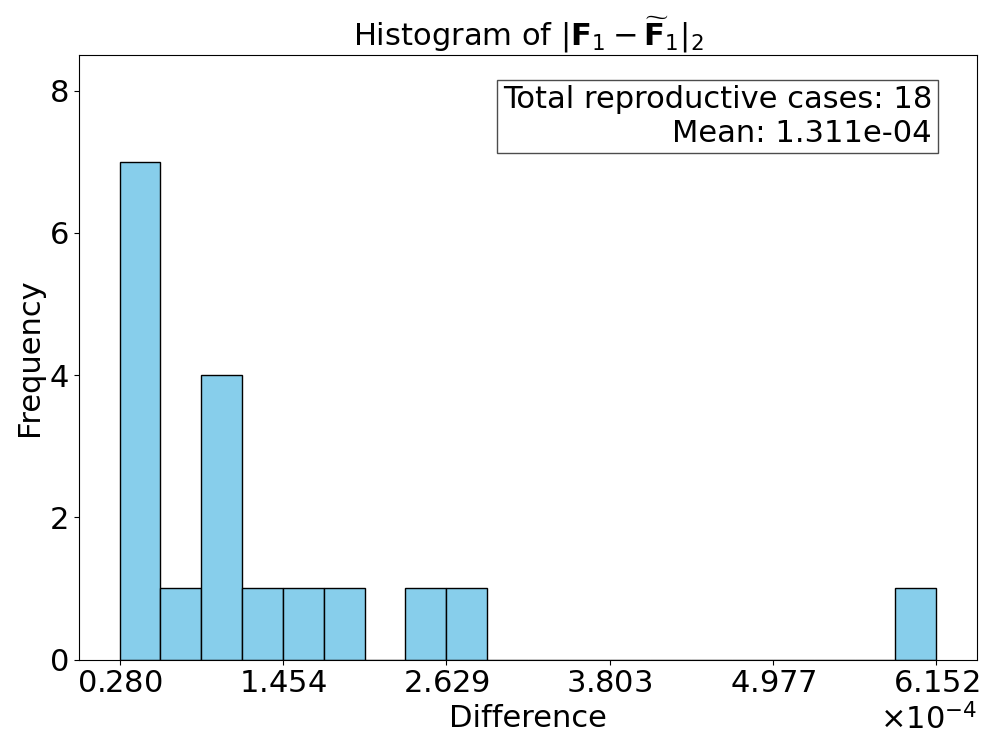}
    \includegraphics[width=0.49\textwidth]{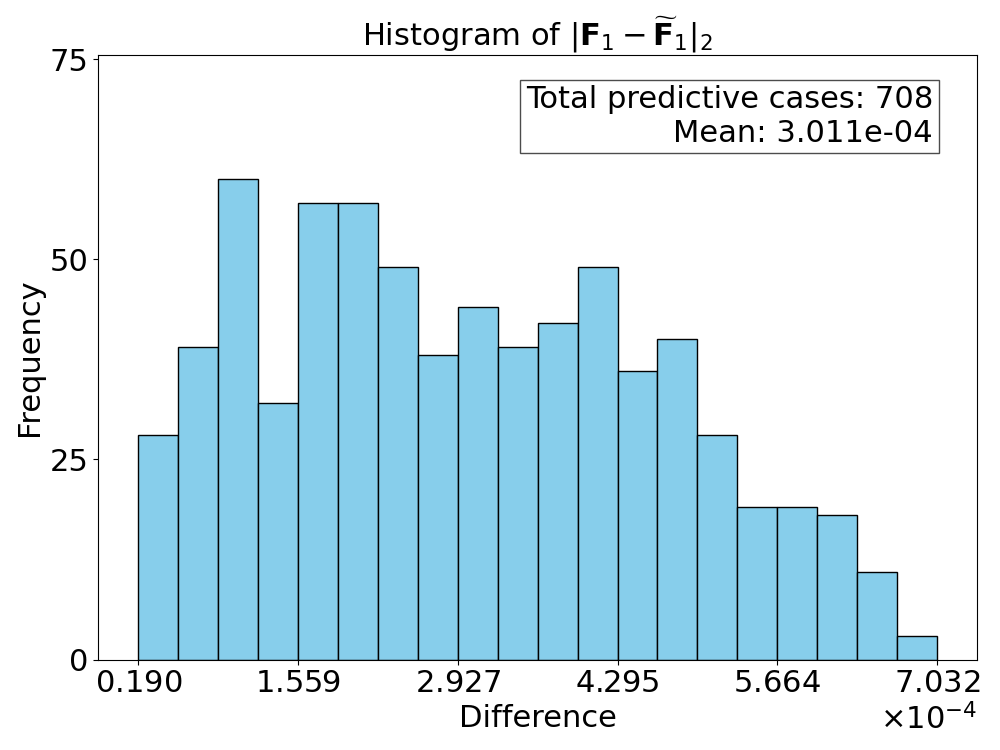} \\
    \includegraphics[width=0.49\textwidth]{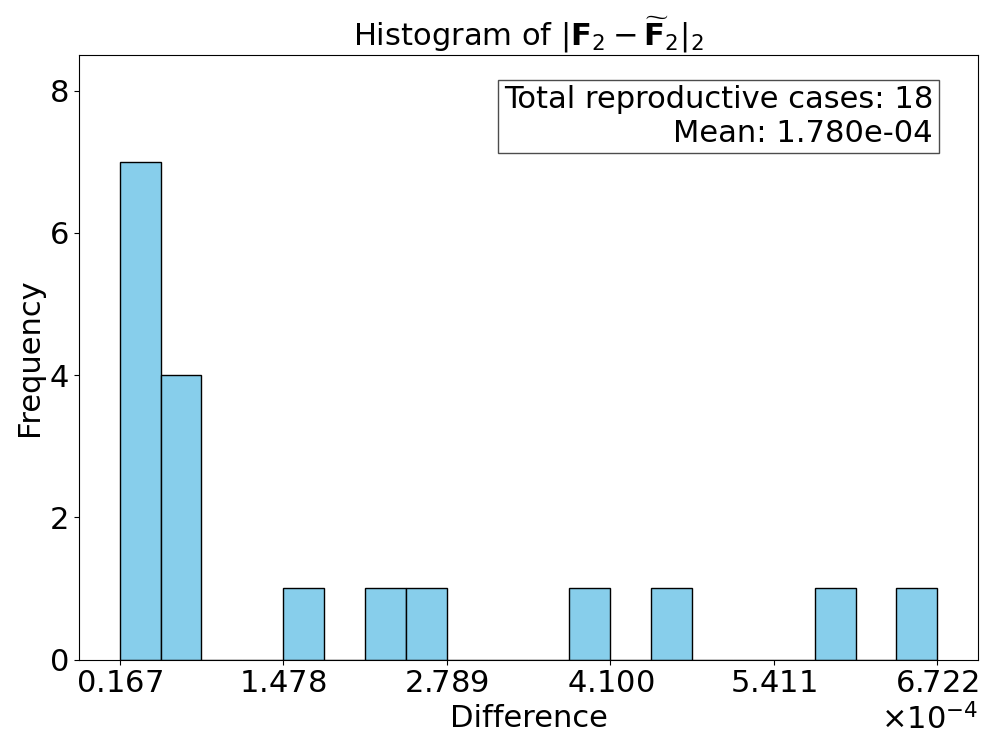}
    \includegraphics[width=0.49\textwidth]{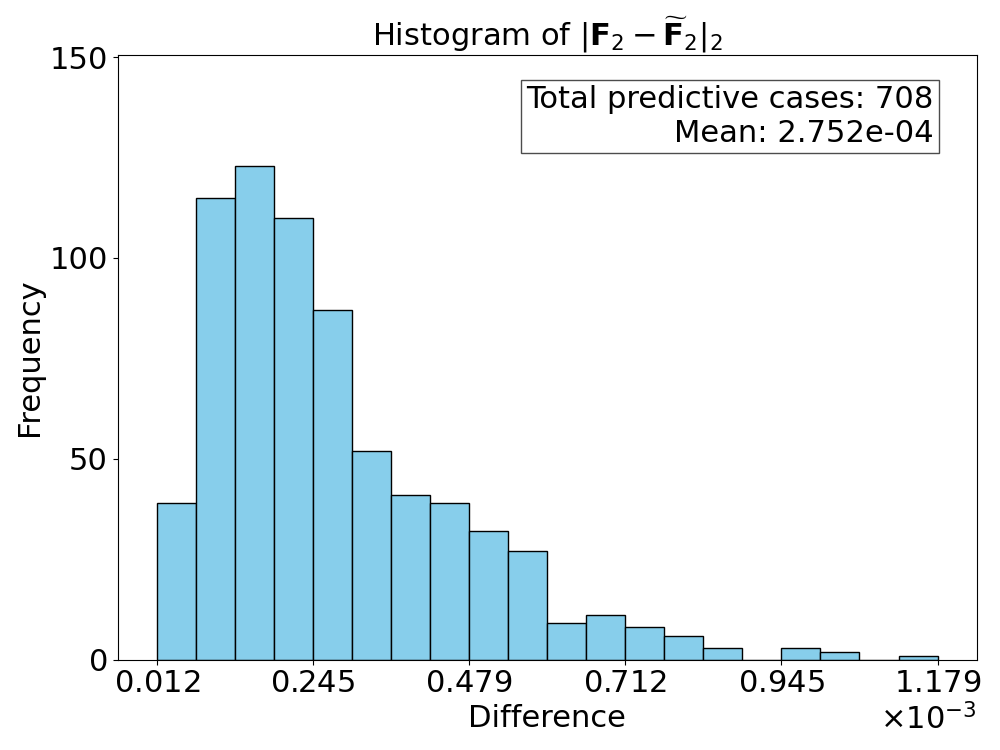}
    \caption{Histogram of the magnitude of the difference in force (in Hartree/Bohr) between FPMD and hybrid AE compression for hydrogen atoms \ce{H1} (top row) and \ce{H2} (bottom row). The left column shows results for the 18 reproductive cases, and the right column shows results for the 708 predictive unseen configurations. The mean force difference is reported in each panel.}
    \label{fig:force_histograms_hybrid_18}
\end{figure}

To assess how compression of the wavefunctions translates into errors in forces, we calculate the approximate forces $\widetilde{\mathbf{F}}$ from the AE-reconstructed wavefunctions for the 18 reproductive configurations and 708 predictive configurations.
Distributions of the corresponding magnitudes of force differences for both hydrogen atoms are shown in Figures \ref{fig:force_histograms_pod_18} and \ref{fig:force_histograms_hybrid_18} for the POD-based and hybrid AEs, respectively, and summarized in Table \ref{tab:force_error_ae}.
The mean force difference is consistently lower for the hybrid AE compared to POD, although the reduction is not as dramatic as that of the error in the snapshot matrix.
While the mean force difference for POD is generally higher than our chosen error tolerance of $5 \times 10^{-4}$ Hartree/Bohr (except for H2 in the predictive configurations, where it is $4.773 \times 10^{-4}$), the mean force difference for the hybrid AE is consistently below this threshold, even for the predictive cases.
However, the maximum force difference for the hybrid AE exceeds the tolerance for both hydrogen atoms in both the reproductive and predictive cases.

\begin{table}[h!]
    \centering
    \begin{tabular}{|c||c|cccc|cccc|}
        \hline
        \multirow{3}{*}{Method} & \multirow{3}{*}{$\frac{\vert \mathbf{Y} - \widetilde{\mathbf{Y}} \vert_\text{F}}{\vert \mathbf{Y} \vert_\text{F}}$} & \multicolumn{4}{c|}{Reproductive} & \multicolumn{4}{c|}{Predictive}  \\
        & & \multicolumn{2}{c}{$\vert \mathbf{F}_1 - \widetilde{\mathbf{F}}_1 \vert_2$} & \multicolumn{2}{c|}{$\vert \mathbf{F}_2 - \widetilde{\mathbf{F}}_2 \vert_2$} & \multicolumn{2}{c}{$\vert \mathbf{F}_1 - \widetilde{\mathbf{F}}_1 \vert_2$} & \multicolumn{2}{c|}{$\vert \mathbf{F}_2 - \widetilde{\mathbf{F}}_2 \vert_2$} \\
        & & Mean & Max. & Mean & Max. & Mean & Max. & Mean & Max. \\
        \hline
        \hline
        POD & 6.462 & 8.042 & 18.160 & 6.632 & 16.512 & 6.119 & 22.151 & 4.773 & 21.096 \\
        Hybrid & 0.518 & 1.311 & 6.152 & 1.780 & 6.722 & 3.011 & 7.032 & 2.752 & 11.788 \\
        %Hybrid & $1.0 \times 10^{-4}$ & $7.2 \times 10^{-5}$ & $5.2 \times 10^{-5}$ & $7.0 \times 10^{-5}$ \\
        \hline
        Factor & 12.470 & 6.135 & 2.952 & 3.726 & 2.456 & 2.032 & 3.150 & 1.735 & 1.790 \\
        \hline
    \end{tabular}
    \caption{
    Relative training error in the snapshot matrix $\mathbf{Y}$ ($\times 10^{-4}$) and mean/maximum magnitude of difference in force (in $\times 10^{-4}$ Hartree/Bohr) from DFT for hydrogen atoms in the reproductive and predictive cases for the POD-based and hybrid AEs. The reduction factor is the ratio of the POD-based AE error to the hybrid AE error.}
    \label{tab:force_error_ae}
\end{table}

Note that this approach to wavefunction compression and force calculation is fundamentally different from our ROM approach, in that we must solve the FOM and then compress/uncompress the wavefunctions rather than directly solving for the wavefunctions in the lower-dimensional ROM subspace.
It is less clear how to efficiently utilize this type of nonlinear compression within a DFT solver, but our results nonetheless illustrate the potential compressive power of nonlinear AE techniques for FPMD applications.

